# Absence of 15-dimensional subalgebras at 16-dimensional Clifford algebra

## Uladzimir Shtukar

Author: Dr. Uladzimir Shtukar, Associate Professor, North Carolina Central University, USA.

e-mail: vshtukar@yahoo.com  Post address: 1906 Raj Drive, Durham, NC 27703.

15A21, 17B30

**Abstract.** It is shown that the $16-$dimensional Clifford algebra over any field has no one $15-$dimensional subalgebra. Canonical bases are used throughout the determination.

Let $V_n$ be the $n-$dimensional vector space with its standard basis $e_1, e_2, \ldots, e_n$. The Clifford algebra over $V_n$ is the associative algebra whose generators are $e_A = e_{i_1} e_{i_2} \ldots e_{i_k}$ for each possible index $A = (i_1 i_2 \ldots i_k)$ that is a subset of $\{1, 2, \ldots, n\}$ with $i_1 < i_2 < \cdots < i_k$. It follows that (as a vector space) the dimension of the Clifford algebra over $V_n$ is $2^n$. The multiplicative identity element of the Clifford algebra is $e_\emptyset$. It is also denoted by $e_0$ and 1; and, it has the following property $e_0 e_i = e_i e_0 = e_i$, $e_0 e_0 = 1$. The vectors $e_1, e_2, \ldots, e_n$ are included in the generators of the Clifford algebra (via index sets $A$ with $k = 1$) and the products of these particular vectors satisfy the following conditions

$$e_i e_j = -e_j e_i, (i \neq j), \text{ and } e_i^2 = -1, \text{ where } i, j = 1, 2, \ldots, n.$$

As mentioned, the multiplication operation of the Clifford algebra is associative, so $(ab)c = a(bc)$ for any elements $a, b, c$ of the Clifford algebra. An arbitrary element of the Clifford algebra over $V_n$ is

$$a = \sum_A a_A e_A,$$

using all indices $A = (i_1 i_2 \ldots i_k)$, where the coefficients $a_A$ are numbers from some field. We denote the Clifford algebra over $V_n$ by $g(n,F)$, when the coefficients $a_A$ are numbers from a given field $F$. For more information about Clifford algebras, see, for example, [1].

Classification of subalgebras is a fundamental way to gain understanding of any algebra. What can we say about subalgebras of the Clifford algebra $g(n,F)$? First, it is easy to see that the Clifford algebra $g(n,F)$ is a $2^n-$dimensional subalgebra of the $2^{n+k}-$dimensional Clifford algebra $g(n+k,F)$ with $k \geq 0$. This leads to the more general and equally obvious statement.

**Lemma. If $h$ is a subalgebra of $2^n-$dimensional Clifford algebra $g(n,F)$, then $h$ is also a subalgebra of $2^{n+k}-$dimensional Clifford algebra $g(n+k,F)$, where $k \geq 0$.**

A lot of properties of Clifford algebras have been established, and some concern subalgebras. For example, equal rank subalgebras are discussed by E. Meinrenken in his book [2]. An answer to the following question would be real progress in the analysis of Clifford algebras: Does the Clifford algebra $g(n,F)$ have $k-$dimensional subalgebras for $2^{n-1} < k < 2^n$? This article answers the question for the case $k = 15$, $n = 4$ by demonstrating the following result:

**Theorem. No one $15-$dimensional subalgebra exists in the $16-$dimensional Clifford algebra $g(4,F)$ over any field $F$.**

**Proof.** For simplicity, the following notation will be used for the generators of $g(4,F)$:

$$e_0 = 1,\ e_1, e_2, e_3, e_4,\ e_{ij} = e_i e_j,\ e_{ijk} = e_i e_j e_k,\ (i,\, j,\, k=1,\, 2,\, 3),\ e_{1234} = e_1 e_2 e_3 e_4.$$

Products of the generators of $g(4,F)$ are recorded in the following tables:

|           | $1$       | $e_1$      | $e_2$        | $e_3$        | $e_4$        | $e_{12}$    | $e_{13}$     | $e_{14}$     |
|-----------|-----------|------------|--------------|--------------|--------------|-------------|--------------|--------------|
| $1$       | $1$       | $e_1$      | $e_2$        | $e_3$        | $e_4$        | $e_{12}$    | $e_{13}$     | $e_{14}$     |
| $e_1$     | $e_1$     | $-1$       | $e_{12}$     | $e_{13}$     | $e_{14}$     | $-e_2$      | $-e_3$       | $-e_4$       |
| $e_2$     | $e_2$     | $-e_{12}$  | $-1$         | $e_{23}$     | $e_{24}$     | $e_1$       | $-e_{123}$   | $-e_{124}$   |
| $e_3$     | $e_3$     | $-e_{13}$  | $-e_{23}$    | $-1$         | $e_{34}$     | $e_{123}$   | $e_1$        | $-e_{134}$   |
| $e_4$     | $e_4$     | $-e_{14}$  | $-e_{24}$    | $-e_{34}$    | $-1$         | $e_{124}$   | $e_{134}$    | $e_1$        |
| $e_{12}$  | $e_{12}$  | $e_2$      | $-e_1$       | $e_{123}$    | $e_{124}$    | $-1$        | $e_{23}$     | $e_{24}$     |
| $e_{13}$  | $e_{13}$  | $e_3$      | $-e_{123}$   | $-e_1$       | $e_{134}$    | $-e_{23}$   | $-1$         | $e_{34}$     |
| $e_{14}$  | $e_{14}$  | $e_4$      | $-e_{124}$   | $-e_{134}$   | $-e_1$       | $-e_{24}$   | $-e_{34}$    | $-1$         |
| $e_{23}$  | $e_{23}$  | $e_{123}$  | $e_3$        | $-e_2$       | $e_{234}$    | $-e_{13}$   | $-e_{12}$    | $e_{1234}$   |
| $e_{24}$  | $e_{24}$  | $e_{124}$  | $e_4$        | $-e_{234}$   | $-e_2$       | $e_{14}$    | $-e_{1234}$  | $-e_{12}$    |
| $e_{34}$  | $e_{34}$  | $e_{134}$  | $e_{234}$    | $e_4$        | $-e_3$       | $e_{1234}$  | $e_{14}$     | $-e_{13}$    |
| $e_{123}$ | $e_{123}$ | $-e_{23}$  | $e_{13}$     | $-e_{12}$    | $e_{1234}$   | $-e_3$      | $e_2$        | $-e_4$       |
| $e_{124}$ | $e_{124}$ | $-e_{24}$  | $e_{14}$     | $-e_{1234}$  | $-e_{12}$    | $-e_4$      | $e_{234}$    | $e_2$        |
| $e_{134}$ | $e_{134}$ | $-e_{34}$  | $e_{1234}$   | $e_{14}$     | $-e_{13}$    | $-e_{234}$  | $-e_4$       | $e_3$        |
| $e_{234}$ | $e_{234}$ | $-e_{1234}$| $-e_{34}$    | $e_{24}$     | $-e_{23}$    | $e_{134}$   | $-e_{124}$   | $e_{123}$    |
| $e_{1234}$| $e_{1234}$| $e_{234}$  | $-e_{134}$   | $e_{124}$    | $-e_{123}$   | $-e_{34}$   | $e_{24}$     | $-e_{23}$    |

|           | $e_{23}$    | $e_{24}$    | $e_{34}$   | $e_{123}$   | $e_{124}$   | $e_{134}$   | $e_{234}$   | $e_{1234}$  |
|-----------|-------------|-------------|------------|-------------|-------------|-------------|-------------|-------------|
| $1$       | $e_{23}$    | $e_{24}$    | $e_{34}$   | $e_{123}$   | $e_{124}$   | $e_{134}$   | $e_{234}$   | $e_{1234}$  |
| $e_1$     | $e_{123}$   | $e_{124}$   | $e_{134}$  | $-e_{23}$   | $-e_{24}$   | $-e_{34}$   | $e_{1234}$  | $-e_{234}$  |
| $e_2$     | $-e_3$      | $-e_4$      | $e_{234}$  | $e_{13}$    | $e_{14}$    | $-e_{1234}$ | $-e_{34}$   | $e_{134}$   |
| $e_3$     | $e_2$       | $-e_{234}$  | $-e_4$     | $-e_{12}$   | $e_{1234}$  | $e_{14}$    | $e_{24}$    | $-e_{124}$  |
| $e_4$     | $e_{234}$   | $e_2$       | $e_3$      | $-e_{1234}$ | $-e_{12}$   | $-e_{13}$   | $-e_{23}$   | $e_{123}$   |
| $e_{12}$  | $-e_{13}$   | $-e_{14}$   | $e_{1234}$ | $-e_3$      | $-e_4$      | $e_{234}$   | $-e_{134}$  | $-e_{34}$   |
| $e_{13}$  | $e_{12}$    | $-e_{1234}$ | $-e_{14}$  | $e_2$       | $-e_{234}$  | $-e_4$      | $e_{124}$   | $e_{24}$    |
| $e_{14}$  | $e_{1234}$  | $e_{12}$    | $e_{13}$   | $e_{234}$   | $e_2$       | $e_3$       | $-e_{123}$  | $-e_{23}$   |
| $e_{23}$  | $-1$        | $e_{34}$    | $-e_{24}$  | $-e_1$      | $e_{34}$    | $-e_{124}$  | $-e_4$      | $-e_{14}$   |
| $e_{24}$  | $-e_{34}$   | $-1$        | $e_{23}$   | $-e_{134}$  | $-e_1$      | $e_{123}$   | $e_3$       | $e_{13}$    |
| $e_{34}$  | $e_{24}$    | $-e_{23}$   | $-1$       | $e_{124}$   | $-e_{123}$  | $-e_1$      | $-e_2$      | $-e_{12}$   |
| $e_{123}$ | $-e_1$      | $e_{134}$   | $-e_{124}$ | $1$         | $-e_{34}$   | $e_{24}$    | $-e_{14}$   | $e_4$       |
| $e_{124}$ | $-e_{134}$  | $-e_1$      | $e_{123}$  | $e_{34}$    | $1$         | $-e_{23}$   | $e_{13}$    | $-e_3$      |
| $e_{134}$ | $e_{124}$   | $-e_{123}$  | $-e_1$     | $-e_{24}$   | $e_{23}$    | $1$         | $-e_{12}$   | $e_2$       |
| $e_{234}$ | $-e_4$      | $e_3$       | $-e_2$     | $e_{14}$    | $-e_{13}$   | $e_{12}$    | $1$         | $-e_1$      |
| $e_{1234}$| $-e_{14}$   | $e_{13}$    | $-e_{12}$  | $-e_4$      | $e_3$       | $-e_2$      | $e_1$       | $1$         |

Our evaluation will be based on canonical bases for $15-$dimensional subspaces of $16-$dimensional vector spaces that were found at [3]. The number of them is 16, and all of them are listed below:

(1) $a_1 = 1 + a_{1,16}e_{1234}$, $a_2 = e_1 + a_{2,16}e_{1234}$, $a_3 = e_2 + a_{3,16}e_{1234}$, $a_4 = e_3 + a_{4,16}e_{1234}$, $a_5 = e_4 + a_{5,16}e_{1234}$, $a_6 = e_{12} + a_{6,16}e_{1234}$, $a_7 = e_{13} + a_{7,16}e_{1234}$, $a_8 = e_{14} + a_{8,16}e_{1234}$, $a_9 = e_{23} + a_{9,16}e_{1234}$, $a_{10} = e_{24} + a_{10,16}e_{1234}$, $a_{11} = e_{34} + a_{11,16}e_{1234}$, $a_{12} = e_{123} + a_{12,16}e_{1234}$, $a_{13} = e_{124} + a_{13,16}e_{1234}$, $a_{14} = e_{134} + a_{14,16}e_{1234}$, $a_{15} = e_{234} + a_{15,16}e_{1234}$.

(2) $a_1 = 1 + a_{1,15}e_{234}$, $a_2 = e_1 + a_{2,15}e_{234}$, $a_3 = e_2 + a_{3,15}e_{234}$, $a_4 = e_3 + a_{4,15}e_{234}$, $a_5 = e_4 + a_{5,15}e_{234}$, $a_6 = e_{12} + a_{6,15}e_{234}$, $a_7 = e_{13} + a_{7,15}e_{234}$, $a_8 = e_{14} + a_{8,15}e_{234}$, $a_9 = e_{23} + a_{9,15}e_{234}$, $a_{10} = e_{24} + a_{10,15}e_{234}$, $a_{11} = e_{34} + a_{11,15}e_{234}$, $a_{12} = e_{123} + a_{12,15}e_{234}$, $a_{13} = e_{124} + a_{13,15}e_{234}$, $a_{14} = e_{134} + a_{14,15}e_{234}$, $a_{15} = e_{1234}$.

(3) $a_1 = 1 + a_{1,14}e_{134}$, $a_2 = e_1 + a_{2,14}e_{134}$, $a_3 = e_2 + a_{3,14}e_{134}$, $a_4 = e_3 + a_{4,14}e_{134}$, $a_5 = e_4 + a_{5,14}e_{134}$, $a_6 = e_{12} + a_{6,14}e_{134}$, $a_7 = e_{13} + a_{7,14}e_{134}$, $a_8 = e_{14} + a_{8,14}e_{134}$, $a_9 = e_{23} + a_{9,14}e_{134}$, $a_{10} = e_{24} + a_{10,14}e_{134}$, $a_{11} = e_{34} + a_{11,14}e_{134}$, $a_{12} = e_{123} + a_{12,14}e_{134}$, $a_{13} = e_{124} + a_{13,14}e_{134}$, $a_{14} = e_{234}$, $a_{15} = e_{1234}$.

(4) $a_1 = 1 + a_{1,13}e_{124}$, $a_2 = e_1 + a_{2,13}e_{124}$, $a_3 = e_2 + a_{3,13}e_{124}$, $a_4 = e_3 + a_{4,13}e_{124}$, $a_5 = e_4 + a_{5,13}e_{124}$, $a_6 = e_{12} + a_{6,13}e_{124}$, $a_7 = e_{13} + a_{7,13}e_{124}$, $a_8 = e_{14} + a_{8,13}e_{124}$, $a_9 = e_{23} + a_{9,13}e_{124}$, $a_{10} = e_{24} + a_{10,13}e_{124}$, $a_{11} = e_{34} + a_{11,13}e_{124}$, $a_{12} = e_{123} + a_{12,13}e_{124}$, $a_{13} = e_{134}$, $a_{14} = e_{234}$, $a_{15} = e_{1234}$.

(5) $a_1 = 1 + a_{1,12}e_{123}$, $a_2 = e_1 + a_{2,12}e_{123}$, $a_3 = e_2 + a_{3,12}e_{123}$, $a_4 = e_3 + a_{4,12}e_{123}$, $a_5 = e_4 + a_{5,12}e_{123}$, $a_6 = e_{12} + a_{6,12}e_{123}$, $a_7 = e_{13} + a_{7,12}e_{123}$, $a_8 = e_{14} + a_{8,12}e_{123}$, $a_9 = e_{23} + a_{9,12}e_{123}$, $a_{10} = e_{24} + a_{10,12}e_{123}$, $a_{11} = e_{34} + a_{11,12}e_{123}$, $a_{12} = e_{124}$, $a_{13} = e_{134}$, $a_{14} = e_{234}$, $a_{15} = e_{1234}$.

(6) $a_1 = 1 + a_{1,11}e_{34}$, $a_2 = e_1 + a_{2,11}e_{34}$, $a_3 = e_2 + a_{3,11}e_{34}$, $a_4 = e_3 + a_{4,11}e_{34}$, $a_5 = e_4 + a_{5,11}e_{34}$, $a_6 = e_{12} + a_{6,11}e_{34}$, $a_7 = e_{13} + a_{7,11}e_{34}$, $a_8 = e_{14} + a_{8,11}e_{34}$, $a_9 = e_{23} + a_{9,11}e_{34}$, $a_{10} = e_{24} + a_{10,11}e_{34}$, $a_{11} = e_{123}$, $a_{12} = e_{124}$, $a_{13} = e_{134}$, $a_{14} = e_{234}$, $a_{15} = e_{1234}$.

(7) $a_1 = 1 + a_{1,10}e_{24}$, $a_2 = e_1 + a_{2,10}e_{24}$, $a_3 = e_2 + a_{3,10}e_{24}$, $a_4 = e_3 + a_{4,10}e_{24}$, $a_5 = e_4 + a_{5,10}e_{24}$, $a_6 = e_{12} + a_{6,10}e_{24}$, $a_7 = e_{13} + a_{7,10}e_{24}$, $a_8 = e_{14} + a_{8,10}e_{24}$, $a_9 = e_{23} + a_{9,10}e_{24}$, $a_{10} = e_{34}$, $a_{11} = e_{123}$, $a_{12} = e_{124}$, $a_{13} = e_{134}$, $a_{14} = e_{234}$, $a_{15} = e_{1234}$.

(8) $a_1 = 1 + a_{1,9}e_{23}$, $a_2 = e_1 + a_{2,9}e_{23}$, $a_3 = e_2 + a_{3,9}e_{23}$, $a_4 = e_3 + a_{4,9}e_{23}$, $a_5 = e_4 + a_{5,9}e_{23}$, $a_6 = e_{12} + a_{6,9}e_{23}$, $a_7 = e_{13} + a_{7,9}e_{23}$, $a_8 = e_{14} + a_{8,9}e_{23}$, $a_9 = e_{24}$, $a_{10} = e_{34}$, $a_{11} = e_{123}$, $a_{12} = e_{124}$, $a_{13} = e_{134}$, $a_{14} = e_{234}$, $a_{15} = e_{1234}$.

(9) $a_1 = 1 + a_{1,8}e_{14}$, $a_2 = e_1 + a_{2,8}e_{14}$, $a_3 = e_2 + a_{3,8}e_{14}$, $a_4 = e_3 + a_{4,8}e_{14}$, $a_5 = e_4 + a_{5,8}e_{14}$, $a_6 = e_{12} + a_{6,8}e_{14}$, $a_7 = e_{13} + a_{7,8}e_{14}$, $a_8 = e_{23}$, $a_9 = e_{24}$, $a_{10} = e_{34}$, $a_{11} = e_{123}$, $a_{12} = e_{124}$, $a_{13} = e_{134}$, $a_{14} = e_{234}$, $a_{15} = e_{1234}$.

(10) $a_1 = 1 + a_{1,7}e_{13}$, $a_2 = e_1 + a_{2,7}e_{13}$, $a_3 = e_2 + a_{3,7}e_{13}$, $a_4 = e_3 + a_{4,7}e_{13}$, $a_5 = e_4 + a_{5,7}e_{13}$, $a_6 = e_{12} + a_{6,7}e_{13}$, $a_7 = e_{14}$, $a_8 = e_{23}$, $a_9 = e_{24}$, $a_{10} = e_{34}$, $a_{11} = e_{123}$, $a_{12} = e_{124}$, $a_{13} = e_{134}$, $a_{14} = e_{234}$, $a_{15} = e_{1234}$.

(11) $a_1 = 1 + a_{1,6}e_{12}$, $a_2 = e_1 + a_{2,6}e_{12}$, $a_3 = e_2 + a_{3,6}e_{12}$, $a_4 = e_3 + a_{4,6}e_{12}$, $a_5 = e_4 + a_{5,6}e_{12}$, $a_6 = e_{13}$, $a_7 = e_{14}$, $a_8 = e_{23}$, $a_9 = e_{24}$, $a_{10} = e_{34}$, $a_{11} = e_{123}$, $a_{12} = e_{124}$, $a_{13} = e_{134}$, $a_{14} = e_{234}$, $a_{15} = e_{1234}$.

(12) $a_1 = 1 + a_{1,5}e_4$, $a_2 = e_1 + a_{2,5}e_4$, $a_3 = e_2 + a_{3,5}e_4$, $a_4 = e_3 + a_{4,5}e_4$, $a_5 = e_{12}$, $a_6 = e_{13}$, $a_7 = e_{14}$, $a_8 = e_{23}$, $a_9 = e_{24}$, $a_{10} = e_{34}$, $a_{11} = e_{123}$, $a_{12} = e_{124}$, $a_{13} = e_{134}$, $a_{14} = e_{234}$, $a_{15} = e_{1234}$.

(13) $a_1 = 1 + a_{1,4}e_3$, $a_2 = e_1 + a_{2,4}e_3$, $a_3 = e_2 + a_{3,4}e_3$, $a_4 = e_4$, $a_5 = e_{12}$, $a_6 = e_{13}$, $a_7 = e_{14}$, $a_8 = e_{23}$, $a_9 = e_{24}$, $a_{10} = e_{34}$, $a_{11} = e_{123}$, $a_{12} = e_{124}$, $a_{13} = e_{134}$, $a_{14} = e_{234}$, $a_{15} = e_{1234}$.

(14)     $a_1 = 1 + a_{1,3}e_2$, $a_2 = e_1 + a_{2,3}e_2$, $a_3 = e_3$, $a_4 = e_4$, $a_5 = e_{12}$, $a_6 = e_{13}$, $a_7 = e_{14}$, $a_8 = e_{23}$, $a_9 = e_{24}$, $a_{10} = e_{34}$, $a_{11} = e_{123}$, $a_{12} = e_{124}$, $a_{13} = e_{134}$, $a_{14} = e_{234}$, $a_{15} = e_{1234}$.

(15)     $a_1 = 1 + a_{1,2}e_1$, $a_2 = e_2$, $a_3 = e_3$, $a_4 = e_4$, $a_5 = e_{12}$, $a_6 = e_{13}$, $a_7 = e_{14}$, $a_8 = e_{23}$, $a_9 = e_{24}$, $a_{10} = e_{34}$, $a_{11} = e_{123}$, $a_{12} = e_{124}$, $a_{13} = e_{134}$, $a_{14} = e_{234}$, $a_{15} = e_{1234}$.

(16)     $a_1 = e_1$, $a_2 = e_2$, $a_3 = e_3$, $a_4 = e_4$, $a_5 = e_{12}$, $a_6 = e_{13}$, $a_7 = e_{14}$, $a_8 = e_{23}$, $a_9 = e_{24}$, $a_{10} = e_{34}$, $a_{11} = e_{123}$, $a_{12} = e_{124}$, $a_{13} = e_{134}$, $a_{14} = e_{234}$, $a_{15} = e_{1234}$.

Every 15-dimensional subalgebra of $g(4,F)$ is a subspace of $g(4,F)$ and hence is associated with exactly one of these 16 canonical bases for some choice of the parameters $a_{ij}$. To determine which subspaces of $g(4,F)$ are subalgebras, we use the well-known fact, that a subspace $h$ of $g(4,F)$ is a subalgebra of $g(4,F)$ if and only if $xy \in h$ for any two elements $x \in h$ and $y \in h$. So, for each of the 16 canonical bases we find all products $a_i a_j$ for $i, j = 1, 2, 3, \ldots, 15$, and then determine whether or not all products $a_i a_j$ are in $h = Span\{a_1, a_2, a_3, \ldots, a_{15}\}$. Specifically, for each of the 16 canonical bases we find all products $a_i a_j$ for $i, j = 1, 2, 3, \ldots, 15$, and then determine whether or not $x_1, x_2, x_3, \ldots, x_{15}$ exist in $F$ so that $a_i a_j = x_1 a_1 + x_2 a_2 + x_3 a_3 + \cdots + x_{15}a_{15}$. The solutions to this equation are found by expressing each $a_i$ in terms of generators of $g(4,F)$ and then equating the coefficients from each side of the equation for corresponding generators. This technique either a) establishes conditions on the parameters $a_i a_j$ for a solution for the $x_i$ to exist or b) establishes that there is no solution. In case a) a subalgebra of $g(4,F)$ is associated with the canonical basis. In case b) no subalgebra of $g(4,F)$ is associated with the canonical basis.

**Basis (1):** $a_1 = 1 + a_{1,16}e_{1234}$, $a_2 = e_1 + a_{2,16}e_{1234}$, $a_3 = e_2 + a_{3,16}e_{1234}$, $a_4 = e_3 + a_{4,16}e_{1234}$, $a_5 = e_4 + a_{5,16}e_{1234}$, $a_6 = e_{12} + a_{6,16}e_{1234}$, $a_7 = e_{13} + a_{7,16}e_{1234}$, $a_8 = e_{14} + a_{8,16}e_{1234}$, $a_9 = e_{23} + a_{9,16}e_{1234}$, $a_{10} = e_{24} + a_{10,16}e_{1234}$, $a_{11} = e_{34} + a_{11,16}e_{1234}$, $a_{12} = e_{123} + a_{12,16}e_{1234}$, $a_{13} = e_{124} + a_{13,16}e_{1234}$, $a_{14} = e_{134} + a_{14,16}e_{1234}$, $a_{15} = e_{234} + a_{15,16}e_{1234}$.

First of all, compute $a_1 a_1$. We have $a_1 a_1 = (1 + a_{1,16}e_{1234})(1 + a_{1,16}e_{1234}) = 1 + 2a_{1,16}e_{1234} + a_{1,16}^2 = x_1 a_1$. So, $x_1 = 1 + a_{1,16}^2$, and $a_{1,16}(a_{1,16}^2 - 1) = 0$. Therefore, $a_{1,16} = 0$ or $a_{1,16} = 1$ or $a_{1,16} = -1$. Consider these 3 cases separately.

Let $a_{1,16} = 0$. Find products $a_2 a_6$, $a_6 a_2$, $a_2 a_{11}$, $a_{11}a_2$, $a_2 a_{13}$, $a_{13}a_2$, $a_3 a_{10}$, $a_{10}a_3$, $a_3 a_{11}$, $a_{11}a_3$, $a_8 a_{10}$, $a_{10}a_8$, $a_3 a_{14}$, $a_{14}a_3$. At this way we will receive a contradiction but it is not the only one way for our goal, other ways with the same result are possible too.

We have $a_2 a_6 = (e_1 + a_{2,16}e_{1234})(e_{12} + a_{6,16}e_{1234}) = -e_2 - a_{6,16}e_{234} - a_{2,16}e_{34} + a_{2,16}a_{6,16} = x_1 a_1 + x_3 a_3 + x_{11}a_{11} + x_{15}a_{15}$. So, $x_1 = a_{2,16}a_{6,16}$, $x_3 = -1$, $x_{11} = -a_{2,16}$, $x_{15} = -a_{6,16}$, and $-a_{3,16} - a_{2,16}a_{11,16} - a_{6,16}a_{15,16} = 0$. Next, $a_6 a_2 = (e_{12} + a_{6,16}e_{1234})(e_1 + a_{2,16}e_{1234}) = e_2 + a_{6,16}e_{234} - a_{2,16}e_{34} + a_{2,16}a_{6,16} = x_1 a_1 + x_3 a_3 + x_{11}a_{11} + x_{15}a_{15}$. So, $x_1 = a_{2,16}a_{6,16}$, $x_3 = 1$, $x_{11} = -a_{2,16}$, $x_{15} = a_{6,16}$, and $a_{3,16} - a_{2,16}a_{11,16} + a_{6,16}a_{15,16} = 0$. Last two equations produce the following results

$$a_{3,16} = -a_{6,16}a_{15,16}, \ a_{2,16}a_{11,16} = 0.$$

We have $a_2 a_{11} = (e_1 + a_{2,16}e_{1234})(e_{34} + a_{11,16}e_{1234}) = e_{134} - a_{11,16}e_{234} - a_{2,16}e_{12} + a_{2,16}a_{11,16} = x_1 a_1 + x_6 a_6 + x_{14}a_{14} + x_{15}a_{15}$. So, $x_1 = a_{2,16}a_{11,16}$, $x_6 = -a_{2,16}$, $x_{14} = 1$,

$x_{15} = -a_{11,16}$, and $-a_{2,16}a_{6,16} + a_{14,16} - a_{11,16}a_{15,16} = 0$. Next, $a_{11}a_2 = (e_{34} + a_{11,16}e_{1234})(e_1 + a_{2,16}e_{1234}) = e_{134} + a_{11,16}e_{234} - a_{2,16}e_{12} + a_{2,16}a_{11,16} = x_1a_1 + x_6a_6 + x_{14}a_{14} + x_{15}a_{15}$. So, $x_1 = a_{2,16}a_{11,16}$, $x_6 = -a_{2,16}$, $x_{14} = 1$, $x_{15} = a_{11,16}$, and $-a_{2,16}a_{6,16} + a_{14,16} + a_{11,16}a_{15,16} = 0$. Last two equations produce the following results

$$a_{14,16} = a_{2,16}a_{6,16}, \ a_{11,16}a_{15,16} = 0.$$

We have $a_2a_{13} = (e_1 + a_{2,16}e_{1234})(e_{124} + a_{13,16}e_{1234}) = -e_{24} - a_{13,16}e_{234} + a_{2,16}e_3 + a_{2,16}a_{13,16} = x_1a_1 + x_4a_4 + x_{10}a_{10} + x_{15}a_{15}$. So, $x_1 = a_{2,16}a_{13,16}$, $x_4 = a_{2,16}$, $x_{10} = -1$, $x_{15} = -a_{13,16}$, and $a_{2,16}a_{4,16} - a_{10,16} - a_{13,16}a_{15,16} = 0$. Next, $a_{13}a_2 = (e_{124} + a_{13,16}e_{1234})(e_1 + a_{2,16}e_{1234}) = -e_{24} + a_{13,16}e_{234} - a_{2,16}e_3 + a_{2,16}a_{13,16} = x_1a_1 + x_4a_4 + x_{10}a_{10} + x_{15}a_{15}$. So, $x_1 = a_{2,16}a_{13,16}$, $x_4 = -a_{2,16}$, $x_{10} = -1$, $x_{15} = a_{13,16}$, and $-a_{2,16}a_{4,16} - a_{10,16} + a_{13,16}a_{15,16} = 0$. Last two equations produce the following results

$$a_{10,16} = 0, \ a_{2,16}a_{4,16} = a_{13,16}a_{15,16}.$$

Compute the product $a_3a_{10} = (e_2 + a_{3,16}e_{1234})(e_{24} + a_{10,16}e_{1234}) = -e_4 + a_{10,16}e_{134} + a_{3,16}e_{13} + a_{3,16}a_{10,16} = x_1a_1 + x_5a_5 + x_7a_7 + x_{14}a_{14}$. Therefore, $x_1 = a_{3,16}a_{10,16}$, $x_5 = -1$, $x_7 = a_{3,16}$, $x_{14} = a_{10,16}$, and $-a_{5,16} + a_{3,16}a_{7,16} + a_{10,16}a_{14,16} = 0$. Next, find $a_{10}a_3 = (e_{24} + a_{10,16}e_{1234})(e_2 + a_{3,16}e_{1234}) = e_4 - a_{10,16}e_{134} + a_{3,16}e_{13} + a_{3,16}a_{10,16} = x_1a_1 + x_5a_5 + x_7a_7 + x_{14}a_{14}$. Therefore, $x_1 = a_{3,16}a_{10,16}$, $x_5 = 1$, $x_7 = a_{3,16}$, $x_{14} = -a_{10,16}$, and $a_{5,16} + a_{3,16}a_{7,16} - a_{10,16}a_{14,16} = 0$. Last two equations generate the following results

$$a_{5,16} = a_{10,16}a_{14,16}, \ a_{3,16}a_{7,16} = 0.$$

We see now that $a_{5,16} = 0$. Compute the product $a_3a_{11} = (e_2 + a_{3,16}e_{1234})(e_{34} + a_{11,16}e_{1234}) = e_{234} + a_{11,16}e_{134} - a_{3,16}e_{12} + a_{3,16}a_{11,16} = x_1a_1 + x_6a_6 + x_{14}a_{14} + x_{15}a_{15}$. Therefore, $x_1 = a_{3,16}a_{11,16}$, $x_6 = -a_{3,16}$, $x_{14} = a_{11,16}$, $x_{15} = 1$, and $a_{15,16} - a_{3,16}a_{6,16} + a_{11,16}a_{14,16} = 0$. Find the next $a_{11}a_3 = (e_{34} + a_{11,16}e_{1234})(e_2 + a_{3,16}e_{1234}) = e_{234} - a_{11,16}e_{134} - a_{3,16}e_{12} + a_{3,16}a_{11,16} = x_1a_1 + x_6a_6 + x_{14}a_{14} + x_{15}a_{15}$. So, $x_1 = a_{3,16}a_{11,16}$, $x_6 = -a_{3,16}$, $x_{14} = -a_{11,16}$, $x_{15} = 1$, and $a_{15,16} - a_{3,16}a_{6,16} - a_{11,16}a_{14,16} = 0$. Last two equations generate the following results

$$a_{15,16} = a_{3,16}a_{6,16}, \ a_{11,16}a_{14,16}{=}0.$$

Compute the product $a_8a_{10} = (e_{14} + a_{8,16}e_{1234})(e_{24} + a_{10,16}e_{1234}) = e_{12} - a_{10,16}e_{23} + a_{8,16}e_{13} + a_{8,16}a_{10,16} = x_1a_1 + x_6a_6 + x_7a_7 + x_9a_9$. Therefore, $x_1 = a_{8,16}a_{10,16}$, $x_6 = 1$, $x_7 = a_{8,16}$, $x_9 = -a_{10,16}$, and $a_{6,16} + a_{7,16}a_{8,16} - a_{10,16}a_{9,16} = 0$. Next, find $a_{10}a_8 = (e_{24} + a_{10,16}e_{1234})(e_{14} + a_{8,16}e_{1234}) = -e_{12} - a_{10,16}e_{23} + a_{8,16}e_{13} + a_{8,16}a_{10,16} = x_1a_1 + x_6a_6 + x_7a_7 + x_9a_9$. Therefore, $x_1 = a_{8,16}a_{10,16}$, $x_6 = -1$, $x_7 = a_{8,16}$, $x_9 = -a_{10,16}$ and $-a_{6,16} + a_{7,16}a_{8,16} - a_{10,16}a_{9,16} = 0$. Last two equations generate the following results

$$a_{6,16} = 0, \ a_{7,16}a_{8,16} = a_{10,16}a_{9,16}.$$

We can see now that $a_{15,16} = 0$, $a_{14,16} = 0$, and $a_{3,16} = 0$ if apply to the previous results.

Compute the product $a_3a_{14} = (e_2 + a_{3,16}e_{1234})(e_{134} + a_{14,16}e_{1234}) = -e_{1234} + a_{14,16}e_{134} - a_{3,16}e_2 + a_{3,16}a_{14,16} = x_1a_1 + x_3a_3 + x_{14}a_{14}$. Therefore, $x_1 = a_{3,16}a_{14,16}$, $x_3 = -a_{3,16}$, $x_{14} = a_{14,16}$, and $-a_{3,16}a_{3,16} + a_{14,16}a_{14,16} = -1$. Product $a_{14}a_3$ generates the same equality.

According the result $a_{3,16} = 0$ that was found before, we obtain $a_{14,16}a_{14,16} = -1$. This property of $a_{14,16}$ contradicts to $a_{14,16} = 0$, and according this fact we have to say that no 15-dimensional subalgebra is generated by Basis (1) in our 16-dimensional Clifford algebra.

Let $a_{1,16} = 1$. Compute products $a_1a_2$, $a_2a_1$. We have $a_1a_2=(1 + e_{1234})(e_1 + a_{2,16}e_{1234}) = e_1 + a_{2,16}e_{1234} + e_{234} + a_{2,16} = x_1a_1 + x_2a_2 + x_{15}a_{15}$. So, $x_1 = a_{2,16}$, $x_2 = 1$, $x_{15} = 1$, and $a_{15,16} = -a_{2,16}$. Similarly, $a_2a_1=(e_1 + a_{2,16}e_{1234})(1 + e_{1234}) = e_1 + a_{2,16}e_{1234} - e_{234} + a_{2,16} = x_1a_1 + x_2a_2 + x_{15}a_{15}$. So, $x_1 = a_{2,16}$, $x_2 = 1$, $x_{15} = -1$, and $a_{15,16} = a_{2,16}$. Last two equations produce $a_{2,16} = 0$, $a_{15,16} = 0$. If we compute now $a_2a_2 = e_1e_1 = -1$ then we see that this product is not located in $Span\{a_1, a_2, a_3, \ldots, a_{16}\}$ because $-1$ can not be written as $x_1a_1$, $x_1a_1 \neq -1$ for any $x_1$. It means that no 15-dimensional subalgebra of 16-dimensional Clifford algebra is generated by Basis (1) at this case.

Let $a_{1,16} = -1$. Compute products $a_1a_2$, $a_2a_1$. We have $a_1a_2=(1 - e_{1234})(e_1 + a_{2,16}e_{1234}) = e_1 + a_{2,16}e_{1234} - e_{234} - a_{2,16} = x_1a_1 + x_2a_2 + x_{15}a_{15}$. So, $x_1 = -a_{2,16}$, $x_2 = 1$, $x_{15} = -1$, and $a_{15,16} = a_{2,16}$. Similarly, $a_2a_1=(e_1 + a_{2,16}e_{1234})(1 - e_{1234}) = e_1 + a_{2,16}e_{1234} + e_{234} + a_{2,16} = x_1a_1 + x_2a_2 + x_{15}a_{15}$. So, $x_1 = -a_{2,16}$, $x_2 = 1$, $x_{15} = 1$, and $a_{15,16} = -a_{2,16}$. Last two equations produce $a_{2,16} = 0$, $a_{15,16} = 0$. If we compute now $a_2a_2 = e_1e_1 = -1$ then we see that this product is not located in $Span\{a_1, a_2, a_3, \ldots, a_{15}\}$ because $-1$ can not be written as $x_1a_1$, $x_1a_1 \neq -1$ for any $x_1$. It means that Basis (1) doesn't generate any 15-dimensional subalgebra of 16-dimensional Clifford algebra at this case.

**Basis (2):** $a_1 = 1 + a_{1,15}e_{234}$, $a_2 = e_1 + a_{2,15}e_{234}$, $a_3 = e_2 + a_{3,15}e_{234}$, $a_4 = e_3 + a_{4,15}e_{234}$, $a_5 = e_4 + a_{5,15}e_{234}$, $a_6 = e_{12} + a_{6,15}e_{234}$, $a_7 = e_{13} + a_{7,15}e_{234}$, $a_8 = e_{14} + a_{8,15}e_{234}$, $a_9 = e_{23} + a_{9,15}e_{234}$, $a_{10} = e_{24} + a_{10,15}e_{234}$, $a_{11} = e_{34} + a_{11,15}e_{234}$, $a_{12} = e_{123} + a_{12,15}e_{234}$, $a_{13} = e_{124} + a_{13,15}e_{234}$, $a_{14} = e_{134} + a_{14,15}e_{234}$, $a_{15} = e_{1234}$

First of all, compute $a_1a_1 = (1 + a_{1,15}e_{234})(1 + a_{1,15}e_{234}) = 1 + a_{1,15}^2 + 2a_{1,15}e_{234} = x_1a_1$. So, $x_1 = 1 + a_{1,15}^2$, and $a_{1,15}(a_{1,15}^2 - 1) = 0$. Therefore, $a_{1,15} = 0$, or $a_{1,15} = 1$, or $a_{1,15} = -1$. Last two cases, $a_{1,15} = 1$ and $a_{1,15} = -1$, are impossible because the product $a_{15}a_{15} = e_{1234}e_{1234} = 1$ can not be written as $x_1a_1$ with these nonzero components $a_{1,15}$. As a fact, we will consider the only case $a_{1,15} = 0$. Start to compute all products $a_ia_j$, $i, j = 2,3, \ldots 15$.

$a_2a_2 = (e_1 + a_{2,15}e_{234})(e_1 + a_{2,15}e_{234}) = -1 + a_{2,15}e_{1234} - a_{2,15}e_{1234} + a_{2,15}^2 = (a_{2,15}^2 - 1)a_1$. It means that $a_2a_2$ belongs to $Span\{a_1, a_2, a_3, \ldots, a_{15}\}$ for any $a_{2,15}$.

$a_2a_3 = (e_1 + a_{2,15}e_{234})(e_2 + a_{3,15}e_{234}) = e_{12} + a_{3,15}e_{1234} - a_{2,15}e_{34} + a_{2,15}a_{3,15} = x_1a_1 + x_6a_6 + x_{11}a_{11} + x_{15}a_{15}$. So, $x_1 = a_{2,15}a_{3,15}$, $x_6 = 1$, $x_{11} = -a_{2,15}$, $x_{15} = a_{3,15}$, and $a_{6,15} - a_{2,15}a_{11,15} = 0$. Vice versa, $a_3a_2 = (e_2 + a_{3,15}e_{234})(e_1 + a_{2,15}e_{234}) = -e_{12} - a_{3,15}e_{1234} - a_{2,15}e_{34} + a_{2,15}a_{3,15} = x_1a_1 + x_6a_6 + x_{11}a_{11}$. So, $x_1 = a_{2,15}a_{3,15}$, $x_6 = -1$, $x_{11} = -a_{2,15}$, $x_{15} = -a_{3,15}$, and $-a_{6,15} - a_{2,15}a_{11,15} = 0$. Last two equations produce the following results $a_{6,15} = 0$, $a_{2,15}a_{11,15} = 0$.

$a_2a_4 = (e_1 + a_{2,15}e_{234})(e_3 + a_{4,15}e_{234}) = e_{13} + a_{4,15}e_{1234} + a_{2,15}e_{24} + a_{2,15}a_{4,15} = x_1a_1 + x_7a_7 + x_{10}a_{10} + x_{15}a_{15}$. So, $x_1 = a_{2,15}a_{4,15}$, $x_7 = 1$, $x_{10} = a_{2,15}$, $x_{15} = a_{4,15}$, and $a_{7,15} + a_{2,15}a_{10,15} = 0$. Vice versa, $a_4a_2 = (e_3 + a_{4,15}e_{234})(e_1 + a_{2,15}e_{234}) = -e_{13} - a_{4,15}e_{1234} + a_{2,15}e_{24} + a_{2,15}a_{4,15} = x_1a_1 + x_7a_7 + x_{10}a_{10} + x_{15}a_{15}$. So, $x_1 = a_{2,15}a_{4,15}$,

$x_7 = -1$, $x_{10} = a_{2,15}$, $x_{15} = -a_{4,15}$, and $-a_{7,15} + a_{2,15}a_{10,15} = 0$. Last two equations produce the following results $a_{7,15} = 0$, $a_{2,15}a_{10,15} = 0$.

$a_2 a_5 = (e_1 + a_{2,15}e_{234})(e_4 + a_{5,15}e_{234}) = e_{14} + a_{5,15}e_{1234} - a_{2,15}e_{23} + a_{2,15}a_{5,15} = x_1 a_1 + x_8 a_8 + x_9 a_9 + x_{15} a_{15}$. So, $x_1 = a_{2,15}a_{5,15}$, $x_8 = 1$, $x_9 = -a_{2,15}$, $x_{15} = a_{5,15}$, and $a_{8,15} - a_{2,15}a_{9,15} = 0$. Vice versa, $a_5 a_2 = (e_4 + a_{5,15}e_{234})(e_1 + a_{2,15}e_{234}) = -e_{14} - a_{5,15}e_{1234} - a_{2,15}e_{24} + a_{2,15}a_{5,15} = x_1 a_1 + x_7 a_7 + x_{10} a_{10} + x_{15} a_{15}$. So, $x_1 = a_{2,15}a_{5,15}$, $x_8 = -1$, $x_9 = -a_{2,15}$, $x_{15} = -a_{5,15}$, and $-a_{8,15} - a_{2,15}a_{9,15} = 0$. Last two equations produce the following results $a_{8,15} = 0$, $a_{2,15}a_{9,15} = 0$.

$a_2 a_6 = (e_1 + a_{2,15}e_{234})(e_{12} + a_{6,15}e_{234}) = -e_2 + a_{6,15}e_{1234} + a_{2,15}e_{134} + a_{2,15}a_{6,15} = x_1 a_1 + x_3 a_3 + x_{14} a_{14} + x_{15} a_{15}$. So, $x_1 = a_{2,15}a_{6,15}$, $x_3 = -1$, $x_{14} = a_{2,15}$, $x_{15} = a_{6,15}$, and $-a_{3,15} + a_{2,15}a_{14,15} = 0$. Vice versa, $a_6 a_2 = (e_{12} + a_{6,15}e_{234})(e_1 + a_{2,15}e_{234}) = e_2 - a_{6,15}e_{1234} - a_{2,15}e_{134} + a_{2,15}a_{6,15} = x_1 a_1 + x_3 a_3 + x_{14} a_{14} + x_{15} a_{15}$. So, $x_1 = a_{2,15}a_{6,15}$, $x_3 = 1$, $x_{11} = -a_{2,15}$, $x_{15} = -a_{6,15}$, and $a_{3,15} - a_{2,15}a_{14,15} = 0$. Last two equations are the same.

$a_2 a_7 = (e_1 + a_{2,15}e_{234})(e_{13} + a_{7,15}e_{234}) = -e_3 + a_{7,15}e_{1234} - a_{2,15}e_{124} + a_{2,15}a_{7,15} = x_1 a_1 + x_4 a_4 + x_{13} a_{13} + x_{15} a_{15}$. So, $x_1 = a_{2,15}a_{7,15}$, $x_4 = -1$, $x_{13} = -a_{2,15}$, $x_{15} = a_{6,15}$, and $-a_{4,15} - a_{2,15}a_{13,15} = 0$. Vice versa, $a_7 a_2 = (e_{13} + a_{7,15}e_{234})(e_1 + a_{2,15}e_{234}) = e_3 - a_{7,15}e_{1234} + a_{2,15}e_{124} + a_{2,15}a_{7,15} = x_1 a_1 + x_4 a_4 + x_{13} a_{13} + x_{15} a_{15}$. So, $x_1 = a_{2,15}a_{7,15}$, $x_4 = 1$, $x_{13} = a_{2,15}$, and $a_{4,15} + a_{2,15}a_{13,15} = 0$. Last two equations are the same.

$a_2 a_8 = (e_1 + a_{2,15}e_{234})(e_{14} + a_{8,15}e_{234}) = -e_4 + a_{8,15}e_{1234} + a_{2,15}e_{123} + a_{2,15}a_{8,15} = x_1 a_1 + x_5 a_5 + x_{12} a_{12} + x_{15} a_{15}$. So, $x_1 = a_{2,15}a_{8,15}$, $x_5 = -1$, $x_{12} = a_{2,15}$, $x_{15} = a_{8,15}$, and $-a_{5,15} + a_{2,15}a_{12,15} = 0$. Product $a_7 a_2$ produces the same equation but details are omitted.

$a_2 a_9 = (e_1 + a_{2,15}e_{234})(e_{23} + a_{9,15}e_{234}) = e_{123} + a_{9,15}e_{1234} - a_{2,15}e_4 + a_{2,15}a_{9,15} = x_1 a_1 + x_5 a_5 + x_{12} a_{12} + x_{15} a_{15}$. So, $x_1 = a_{2,15}a_{9,15}$, $x_5 = -a_{2,15}$, $x_{12} = 1$, $x_{15} = a_{9,15}$, and $-a_{2,15}a_{5,15} + a_{12,15} = 0$. Product $a_9 a_2$ produces the same equation but details are omitted.

$a_2 a_{10} = (e_1 + a_{2,15}e_{234})(e_{24} + a_{10,15}e_{234}) = e_{124} + a_{10,15}e_{1234} + a_{2,15}e_3 + a_{2,15}a_{10,15} = x_1 a_1 + x_4 a_4 + x_{13} a_{13} + x_{15} a_{15}$. So, $x_1 = a_{2,15}a_{10,15}$, $x_4 = a_{2,15}$, $x_{13} = 1$, $x_{15} = a_{10,15}$, and $a_{2,15}a_{4,15} + a_{13,15} = 0$. Product $a_{10} a_2$ produces the same equation but details are omitted.

$a_2 a_{11} = (e_1 + a_{2,15}e_{234})(e_{34} + a_{11,15}e_{234}) = e_{134} + a_{11,15}e_{1234} - a_{2,15}e_2 + a_{2,15}a_{11,15} = x_1 a_1 + x_3 a_3 + x_{14} a_{14} + x_{15} a_{15}$. So, $x_1 = a_{2,15}a_{11,15}$, $x_3 = -a_{2,15}$, $x_{14} = 1$, $x_{15} = a_{11,15}$, and $-a_{2,15}a_{3,15} + a_{14,15} = 0$. Product $a_{11} a_2$ produces the same equation but details are omitted.

$a_2 a_{12} = (e_1 + a_{2,15}e_{234})(e_{123} + a_{12,15}e_{234}) = -e_{23} + a_{12,15}e_{1234} + a_{2,15}e_{14} + a_{2,15}a_{12,15} = x_1 a_1 + x_8 a_8 + x_9 a_9 + x_{15} a_{15}$. So, $x_1 = a_{2,15}a_{12,15}$, $x_8 = a_{2,15}$, $x_9 = -1$, $x_{15} = a_{12,15}$, and $a_{2,15}a_{8,15} - a_{9,15} = 0$. Compute now $a_{12} a_2 = (e_{123} + a_{12,15}e_{234})(e_1 + a_{2,15}e_{234}) = -e_{23} - a_{12,15}e_{1234} - a_{2,15}e_{14} + a_{2,15}a_{12,15} = x_1 a_1 + x_8 a_8 + x_9 a_9 + x_{15} a_{15}$. So, $x_1 = a_{2,15}a_{12,15}$, $x_8 = -a_{2,15}$, $x_9 = -1$, $x_{15} = -a_{12,15}$, and $-a_{2,15}a_{8,15} - a_{9,15} = 0$. Last two equations produce the following results $a_{9,15} = 0$, $a_{2,15}a_{8,15} = 0$.

$a_2 a_{13} = (e_1 + a_{2,15}e_{234})(e_{124} + a_{13,15}e_{234}) = -e_{24} + a_{13,15}e_{1234} - a_{2,15}e_{13} + a_{2,15}a_{13,15} = x_1 a_1 + x_7 a_7 + x_{10}a_{10} + x_{15}a_{15}$. So, $x_1 = a_{2,15}a_{13,15}$, $x_7 = -a_{2,15}$, $x_{10} = -1$, $x_{15} = a_{13,15}$, and $-a_{2,15}a_{7,15} - a_{10,15} = 0$. Compute now $a_{13}a_2 = (e_{124} + a_{13,15}e_{234})(e_1 + a_{2,15}e_{234}) = -e_{24} - a_{13,15}e_{1234} + a_{2,15}e_{13} + a_{2,15}a_{13,15} = x_1 a_1 + x_7 a_7 + x_{10}a_{10} + x_{15}a_{15}$. So, $x_1 = a_{2,15}a_{13,15}$, $x_7 = a_{2,15}$, $x_{10} = -1$, $x_{15} = -a_{13,15}$, and $a_{2,15}a_{7,15} - a_{10,15} = 0$. Last two equations produce the following results $a_{10,15} = 0$, $a_{2,15}a_{7,15} = 0$.

$a_2 a_{14} = (e_1 + a_{2,15}e_{234})(e_{134} + a_{14,15}e_{234}) = -e_{34} + a_{14,15}e_{1234} + a_{2,15}e_{12} + a_{2,15}a_{14,15} = x_1 a_1 + x_6 a_6 + x_{11}a_{11} + x_{15}a_{15}$. So, $x_1 = a_{2,15}a_{14,15}$, $x_6 = a_{2,15}$, $x_{11} = -1$, $x_{15} = a_{14,15}$, and $a_{2,15}a_{6,15} - a_{11,15} = 0$. Compute now $a_{14}a_2 = (e_{134} + a_{14,15}e_{234})(e_1 + a_{2,15}e_{234}) = -e_{34} - a_{14,15}e_{1234} - a_{2,15}e_{12} + a_{2,15}a_{14,15} = x_1 a_1 + x_6 a_6 + x_{11}a_{11} + x_{15}a_{15}$. So, $x_1 = a_{2,15}a_{14,15}$, $x_6 = -a_{2,15}$, $x_{11} = -1$, $x_{15} = -a_{14,15}$, and $-a_{2,15}a_{6,15} - a_{11,15} = 0$. Last two equations produce the following results $a_{11,15} = 0$, $a_{2,15}a_{6,15} = 0$.

$a_2 a_{15} = (e_1 + a_{2,15}e_{234})(e_{1234}) = -e_{234} - a_{2,15}e_1 = x_2 a_2$. So, $x_2 = -a_{2,15}$, and $a_{2,15}^2 = 1$. Finally, $a_{2,15} = 1$ or $a_{2,15} = -1$.

According our results for products $a_2 a_6$, $a_6 a_2$, we can determine that $a_{3,15} = 0$. The following two products $a_3 a_6$ and $a_6 a_3$ play the crucial role for Basis (2). Find them.

$a_3 a_6 = (e_2 + a_{3,15}e_{234})(e_{12} + a_{6,15}e_{234}) = e_1 - a_{6,15}e_{34} + a_{3,15}e_{134} + a_{3,15}a_{6,15} = x_1 a_1 + x_2 a_2 + x_{11}a_{11} + x_{14}a_{14}$. So, $x_1 = a_{3,15}a_{6,15}$, $x_2 = 1$, $x_{11} = -a_{6,15}$, $x_{14} = a_{3,15}$, and $a_{2,15} - a_{6,15}a_{11,15} + a_{3,15}a_{14,15} = 0$. Similarly, $a_6 a_3 = (e_{12} + a_{6,15}e_{234})(e_2 + a_{3,15}e_{234}) = -e_1 - a_{6,15}e_{34} - a_{3,15}e_{134} + a_{3,15}a_{6,15} = x_1 a_1 + x_2 a_2 + x_{11}a_{11} + x_{14}a_{14}$. So, $x_1 = a_{3,15}a_{6,15}$, $x_2 = -1$, $x_{11} = -a_{6,15}$, $x_{14} = -a_{3,15}$, and $-a_{2,15} - a_{6,15}a_{11,15} - a_{3,15}a_{14,15} = 0$. Last two equations produce the following results $a_{6,15}a_{11,15} = 0$, $a_{2,15} + a_{3,15}a_{14,15} = 0$, and therefore $a_{2,15} = 0$ which contradicts to $a_{2,15} = 1$ or $a_{2,15} = -1$. This contradiction says that Basis (2) doesn't generate any 15-dimensional subalgebra in the 16-dimensional Clifford algebra.

**Basis (3):** $a_1 = 1 + a_{1,14}e_{134}$, $a_2 = e_1 + a_{2,14}e_{134}$, $a_3 = e_2 + a_{3,14}e_{134}$, $a_4 = e_3 + a_{4,14}e_{134}$, $a_5 = e_4 + a_{5,14}e_{134}$, $a_6 = e_{12} + a_{6,14}e_{134}$, $a_7 = e_{13} + a_{7,14}e_{134}$, $a_8 = e_{14} + a_{8,14}e_{134}$, $a_9 = e_{23} + a_{9,14}e_{134}$, $a_{10} = e_{24} + a_{10,14}e_{134}$, $a_{11} = e_{34} + a_{11,14}e_{134}$, $a_{12} = e_{123} + a_{12,14}e_{134}$, $a_{13} = e_{124} + a_{13,14}e_{134}$, $a_{14} = e_{234}$, $a_{15} = e_{1234}$.

First of all, compute $a_1 a_1 = (1 + a_{1,14}e_{134})(1 + a_{1,14}e_{134}) = 1 + a_{1,14}^2 + 2a_{1,14}e_{134} = x_1 a_1$. So, $x_1 = 1 + a_{1,14}^2$, and $a_{1,14}(a_{1,14}^2 - 1) = 0$. Therefore, $a_{1,14} = 0$, or $a_{1,14} = 1$, or $a_{1,14} = -1$. Last two cases, $a_{1,14} = 1$ and $a_{1,14} = -1$, are impossible because the product $a_{15}a_{15} = e_{1234}e_{1234} = 1$ can not be written as $x_1 a_1$ with these nonzero components $a_{1,14}$. As a fact, we will consider the only case $a_{1,14} = 0$ in our evaluation below. Start to compute products $a_i a_j$, $i, j = 2,3,\ldots 15$.

$a_2 a_2 = (e_1 + a_{2,14}e_{134})(e_1 + a_{2,14}e_{134}) = -1 - a_{2,14}e_{34} - a_{2,14}e_{34} + a_{2,14}^2 = x_1 a_1 + x_{11}a_{11}$. So, $x_1 = a_{2,14}^2 - 1$, $x_{11} = -2a_{2,14}$, and $2a_{2,14}a_{11,14} = 0$. It means that $a_{2,14} = 0$ or $a_{11,14} = 0$.

$a_2 a_3 = (e_1 + a_{2,14}e_{134})(e_2 + a_{3,14}e_{134}) = e_{12} - a_{3,14}e_{34} + a_{2,14}e_{1234} + a_{2,14}a_{3,14} = x_1 a_1 + x_6 a_6 + x_{11}a_{11} + x_{15}a_{15}$. So, $x_1 = a_{2,14}a_{3,14}$, $x_6 = 1$, $x_{11} = -a_{3,14}$, $x_{15} = a_{2,14}$, and $a_{6,14} - a_{3,14}a_{11,14} = 0$. Vice versa, $a_3 a_2 = (e_2 + a_{3,14}e_{134})(e_1 + a_{2,14}e_{134}) = -e_{12} -$

$a_{2,14}e_{1234} - a_{3,14}e_{34} + a_{2,14}a_{3,14} = x_1a_1 + x_6a_6 + x_{11}a_{11} + x_{15}a_{15}$. So, $x_1 = a_{2,14}a_{3,14}$, $x_6 = -1$, $x_{11} = -a_{3,14}$, $x_{15} = -a_{2,14}$, and $-a_{6,14} - a_{3,14}a_{11,14} = 0$. Last two equations produce the following results $a_{6,14} = 0$, $a_{2,14}a_{11,14} = 0$.

$a_2a_4 = (e_1 + a_{2,14}e_{134})(e_3 + a_{4,14}e_{134}) = e_{13} - a_{4,14}e_{34} + a_{2,14}e_{14} + a_{2,14}a_{4,14} = x_1a_1 + x_7a_7 + x_8a_8 + x_{11}a_{11}$. So, $x_1 = a_{2,14}a_{4,14}$, $x_7 = 1$, $x_8 = a_{2,14}$, $x_{11} = -a_{4,14}$, and $a_{7,14} + a_{2,14}a_{8,14} - a_{4,14}a_{11,14} = 0$. Vice versa, $a_4a_2 = (e_3 + a_{4,14}e_{134})(e_1 + a_{2,14}e_{134}) = -e_{13} - a_{4,14}e_{34} + a_{2,14}e_{14} + a_{2,14}a_{4,14} = x_1a_1 + x_7a_7 + x_8a_8 + x_{11}a_{11}$. So, $x_1 = a_{2,14}a_{4,14}$, $x_7 = -1$, $x_8 = a_{2,14}$, $x_{11} = -a_{4,14}$, and $-a_{7,14} + a_{2,14}a_{8,14} - a_{4,14}a_{11,14} = 0$. Last two equations produce the following results $a_{7,14} = 0$, $a_{2,14}a_{8,14} = a_{4,14}a_{11,14}$.

$a_2a_5 = (e_1 + a_{2,14}e_{134})(e_4 + a_{5,14}e_{134}) = e_{14} - a_{5,14}e_{34} - a_{2,14}e_{13} + a_{2,14}a_{5,14} = x_1a_1 + x_7a_7 + x_8a_8 + x_{11}a_{11}$. So, $x_1 = a_{2,14}a_{5,14}$, $x_7 = -a_{2,14}$, $x_8 = 1$, $x_{11} = -a_{5,14}$, and $-a_{2,14}a_{7,14} + a_{8,14} - a_{5,14}a_{11,14} = 0$. Vice versa, $a_5a_2 = (e_4 + a_{5,14}e_{134})(e_1 + a_{2,14}e_{134}) = -e_{14} - a_{5,14}e_{34} - a_{2,14}e_{13} + a_{2,14}a_{5,14} = x_1a_1 + x_7a_7 + x_8a_8 + x_{11}a_{11}$. So, $x_1 = a_{2,14}a_{5,14}$, $x_7 = -a_{2,14}$, $x_8 = -1$, $x_{11} = -a_{5,14}$, and $-a_{2,14}a_{7,14} - a_{8,14} - a_{5,14}a_{11,14} = 0$.. Last two equations produce the following results $a_{8,14} = 0$, $a_{2,14}a_{7,14} = -a_{5,14}a_{11,14}$.

$a_2a_6 = (e_1 + a_{2,14}e_{134})(e_{12} + a_{6,14}e_{134}) = -e_2 - a_{6,14}e_{34} - a_{2,14}e_{234} + a_{2,14}a_{6,14} = x_1a_1 + x_3a_3 + x_{11}a_{11} + x_{14}a_{14}$. So, $x_1 = a_{2,14}a_{6,14}$, $x_3 = -1$, $x_{11} = -a_{6,14}$, $x_{14} = -a_{2,14}$, and $-a_{3,14} - a_{6,14}a_{11,14} = 0$. Vice versa, $a_6a_2 = (e_{12} + a_{6,14}e_{134})(e_1 + a_{2,14}e_{134}) = e_2 - a_{6,14}e_{34} + a_{2,14}e_{234} + a_{2,14}a_{6,14} = x_1a_1 + x_3a_3 + x_{11}a_{11} + x_{14}a_{14}$. So, $x_1 = a_{2,14}a_{6,14}$, $x_3 = 1$, $x_{11} = -a_{6,14}$, $x_{14} = a_{2,14}$, and $a_{3,14} - a_{6,14}a_{11,14} = 0$. Last two equations produce the following results $a_{3,14} = 0$, $a_{6,14}a_{11,14} = 0$.

$a_2a_7 = (e_1 + a_{2,14}e_{134})(e_{13} + a_{7,14}e_{134}) = -e_3 - a_{7,14}e_{34} - a_{2,14}e_4 + a_{2,14}a_{7,14} = x_1a_1 + x_4a_4 + x_5a_5 + x_{11}a_{11}$. So, $x_1 = a_{2,14}a_{7,14}$, $x_4 = -1$, $x_5 = -a_{2,14}$, $x_{11} = -a_{7,14}$, and $-a_{4,14} - a_{2,14}a_{5,14} - a_{7,14}a_{11,14} = 0$. Vice versa, $a_7a_2 = (e_{13} + a_{7,14}e_{134})(e_1 + a_{2,14}e_{134}) = e_3 - a_{7,14}e_{34} - a_{2,14}e_4 + a_{2,14}a_{7,14} = x_1a_1 + x_4a_4 + x_5a_5 + x_{11}a_{11}$. So, $x_1 = a_{2,14}a_{7,14}$, $x_4 = 1$, $x_5 = -a_{2,14}$, $x_{11} = -a_{7,14}$, and $a_{4,14} - a_{2,14}a_{5,14} - a_{7,14}a_{11,14} = 0$. Last two equations produce the following results $a_{4,14} = 0$, $a_{2,14}a_{5,14} = -a_{7,14}a_{11,14}$.

$a_2a_8 = (e_1 + a_{2,14}e_{134})(e_{14} + a_{8,14}e_{134}) = -e_4 - a_{8,14}e_{34} + a_{2,14}e_3 + a_{2,14}a_{8,14} = x_1a_1 + x_4a_4 + x_5a_5 + x_{11}a_{11}$. So, $x_1 = a_{2,14}a_{8,14}$, $x_4 = a_{2,14}$, $x_5 = -1$, $x_{11} = -a_{8,14}$, and $a_{2,14}a_{4,14} - a_{5,14} - a_{8,14}a_{11,14} = 0$. Product $a_8a_2$ gives $a_8a_2 = (e_{14} + a_{8,14}e_{134})(e_1 + a_{2,14}e_{134}) = e_4 + a_{2,14}e_3 - a_{8,14}e_{34} + a_{2,14}a_{8,14} = x_1a_1 + x_4a_4 + x_5a_5 + x_{11}a_{11}$. So, $x_1 = a_{2,14}a_{8,14}$, $x_4 = a_{2,14}$, $x_5 = 1$, $x_{11} = -a_{8,14}$, and $a_{2,14}a_{4,14} + a_{5,14} - a_{8,14}a_{11,14} = 0$. Last two equations produce the following results $a_{5,14} = 0$, $a_{2,14}a_{5,14} = a_{8,14}a_{11,14}$.

$a_2a_9 = (e_1 + a_{2,14}e_{134})(e_{23} + a_{9,14}e_{134}) = e_{123} - a_{9,14}e_{34} + a_{2,14}e_{124} + a_{2,14}a_{9,14} = x_1a_1 + x_{11}a_{11} + x_{12}a_{12} + x_{13}a_{13}$. So, $x_1 = a_{2,14}a_{9,14}$, $x_{11} = -a_{9,14}$, $x_{12} = 1$, $x_{13} = a_{2,14}$, and $-a_{9,14}a_{11,14} + a_{12,14} + a_{2,14}a_{13,14} = 0$. Vice versa, $a_9a_2 = (e_{23} + a_{9,14}e_{134})(e_1 + a_{2,14}e_{134}) = e_{123} - a_{9,14}e_{34} - a_{2,14}e_{124} + a_{2,14}a_{9,14} = x_1a_1 + x_{11}a_{11} + x_{12}a_{12} + x_{13}a_{13}$. So, $x_1 = a_{2,14}a_{9,14}$, $x_{11} = -a_{9,14}$, $x_{12} = 1$, $x_{13} = -a_{2,14}$, and $-a_{9,14}a_{11,14} + a_{12,14} - a_{2,14}a_{13,14} = 0$. Last two equations produce the following results $a_{2,14}a_{13,14} = 0$, $a_{12,14} = a_{9,14}a_{11,14}$.

$a_2a_{10} = (e_1 + a_{2,14}e_{134})(e_{24} + a_{10,14}e_{134}) = e_{124} - a_{10,14}e_{34} - a_{2,14}e_{123} + a_{2,14}a_{10,14} = x_1a_1 + x_{11}a_{11} + x_{12}a_{12} + x_{13}a_{13}$. So, $x_1 = a_{2,14}a_{10,14}$, $x_{11} = -a_{10,14}$, $x_{12} = -a_{2,14}$, $x_{13} = 1$, and $-a_{10,14}a_{11,14} - a_{2,14}a_{12,14} + a_{13,14} = 0$. Vice versa, $a_{10}a_2 = (e_{24} + a_{10,14}e_{134})(e_1 + a_{2,14}e_{134}) = e_{124} - a_{10,14}e_{34} + a_{2,14}e_{123} + a_{2,14}a_{10,14} = x_1a_1 + x_{11}a_{11} + x_{12}a_{12} + x_{13}a_{13}$. So, $x_1 = a_{2,14}a_{10,14}$, $x_{11} = -a_{10,14}$, $x_{12} = a_{2,14}$, $x_{13} = 1$, and $-a_{10,14}a_{11,14} + a_{2,14}a_{12,14} + a_{13,14} = 0$. Last two equations produce the following results $a_{2,14}a_{12,14} = 0$, $a_{13,14} = a_{10,14}a_{11,14}$.

$a_2a_{11} = (e_1 + a_{2,14}e_{134})(e_{34} + a_{11,14}e_{134}) = e_{134} - a_{11,14}e_{34} - a_{2,14}e_1 + a_{2,14}a_{11,14} = x_1a_1 + x_2a_2 + x_{11}a_{11} + x_{14}a_{14}$. So, $x_1 = a_{2,14}a_{11,14}$, $x_2 = -a_{2,14}$, $x_{11} = -a_{11,14}$, and $-a_{2,14}a_{2,14} - a_{11,14}a_{11,14} = 1$. Vice versa, $a_{11}a_2 = (e_{34} + a_{11,14}e_{134})(e_1 + a_{2,14}e_{134}) = e_{134} - a_{11,14}e_{34} - a_{2,14}e_1 + a_{2,14}a_{11,14} = a_2a_{11}$. So, the product $a_{11}a_2$ produces the same equation as $a_2a_{11}$ does.

$a_2a_{12} = (e_1 + a_{2,14}e_{134})(e_{123} + a_{12,14}e_{134}) = -e_{23} - a_{12,14}e_{34} - a_{2,14}e_{24} + a_{2,14}a_{12,14} = x_1a_1 + x_9a_9 + x_{10}a_{10} + x_{11}a_{11}$. So, $x_1 = a_{2,14}a_{12,14}$, $x_9 = -1$, $x_{10} = -a_{2,14}$, $x_{11} = -a_{12,14}$, and $-a_{9,14} - a_{2,14}a_{10,14} - a_{11,14}a_{12,14} = 0$. Compute $a_{12}a_2 = (e_{123} + a_{12,14}e_{134})(e_1 + a_{2,14}e_{134}) = -e_{23} - a_{12,14}e_{34} + a_{2,14}e_{24} + a_{2,14}a_{12,14} = x_1a_1 + x_9a_9 + x_{10}a_{10} + x_{11}a_{11}$. So, $x_1 = a_{2,14}a_{12,14}$, $x_9 = -1$, $x_{10} = a_{2,14}$, $x_{11} = -a_{12,14}$, and $-a_{9,14} + a_{2,14}a_{10,14} - a_{11,14}a_{12,14} = 0$. Last two equations produce the following results $a_{9,14} = -a_{11,14}a_{12,14}$, $a_{2,14}a_{10,14} = 0$.

$a_2a_{13} = (e_1 + a_{2,14}e_{134})(e_{124} + a_{13,14}e_{134}) = -e_{24} - a_{13,14}e_{34} + a_{2,14}e_{23} + a_{2,14}a_{13,14} = x_1a_1 + x_9a_9 + x_{10}a_{10} + x_{11}a_{11}$. So, $x_1 = a_{2,14}a_{13,14}$, $x_9 = a_{2,14}$, $x_{10} = -1$, $x_{11} = -a_{13,14}$, and $a_{2,14}a_{9,14} - a_{10,14} - a_{13,14}a_{11,14} = 0$. Compute now $a_{13}a_2 = (e_{124} + a_{13,14}e_{134})(e_1 + a_{2,14}e_{134}) = -e_{24} - a_{13,14}e_{34} - a_{2,14}e_{23} + a_{2,14}a_{13,14} = x_1a_1 + x_9a_9 + x_{10}a_{10} + x_{11}a_{11}$. So, $x_1 = a_{2,14}a_{13,14}$, $x_9 = -a_{2,14}$, $x_{10} = -1$, $x_{11} = -a_{13,14}$, and $-a_{2,14}a_{9,14} - a_{10,14} - a_{13,14}a_{11,14} = 0$. Last two equations produce the following results $a_{10,14} = -a_{13,14}a_{11,14}$, $a_{2,14}a_{9,14} = 0$.

$a_2a_{14} = (e_1 + a_{2,14}e_{134})(e_{234}) = e_{1234} - a_{2,14}e_{12} = x_6a_6 + x_{15}a_{15}$. So, $x_6 = -a_{2,14}$, $x_{15} = 1$, and $-a_{2,14}a_{6,14} = 0$. Product $a_{14}a_2$ generates the same equation.

$a_2a_{15} = (e_1 + a_{2,14}e_{134})(e_{1234}) = -e_{234} + a_{2,14}e_2 = x_3a_3 + x_{14}a_{14}$. So, $x_3 = a_{2,14}$, $x_{14} = -1$, and $a_{2,14}a_{3,14} = 0$.

Our evaluation shows that $a_{3,14} = 0$, $a_{4,14} = 0$, $a_{5,14} = 0$, $a_{6,14} = 0$, $a_{7,14} = 0$, $a_{8,14} = 0$. Compute now $a_4a_8 = e_3e_{14} = -e_{134}$. But this element $-e_{134}$ can not be written as a linear combination of vectors from the Basis (3). It means that basis (3) doesn't generate any 15-dimensional subalgebra in 16-dimensional Clifford algebra.

**Basis (4):** $a_1 = 1 + a_{1,13}e_{124}$, $a_2 = e_1 + a_{2,13}e_{124}$, $a_3 = e_2 + a_{3,13}e_{124}$, $a_4 = e_3 + a_{4,13}e_{124}$, $a_5 = e_4 + a_{5,13}e_{124}$, $a_6 = e_{12} + a_{6,13}e_{124}$, $a_7 = e_{13} + a_{7,13}e_{124}$, $a_8 = e_{14} + a_{8,13}e_{124}$, $a_9 = e_{23} + a_{9,13}e_{124}$, $a_{10} = e_{24} + a_{10,13}e_{124}$, $a_{11} = e_{34} + a_{11,13}e_{124}$, $a_{12} = e_{123} + a_{12,13}e_{124}$, $a_{13} = e_{134}$, $a_{14} = e_{234}$, $a_{15} = e_{1234}$.

First of all, compute $a_1a_1 = (1 + a_{1,13}e_{124})(1 + a_{1,13}e_{124}) = 1 + a_{1,13}^2 + 2a_{1,13}e_{124} = x_1a_1$. So, $x_1 = 1 + a_{1,13}^2$, and $a_{1,13}(a_{1,13}^2 - 1) = 0$. Therefore, $a_{1,13} = 0$, or $a_{1,13} = 1$, or $a_{1,13} =$ $-1$. Last two cases, $a_{1,13} = 1$ and $a_{1,13} = -1$, are impossible because the product $a_{15}a_{15} = e_{1234}e_{1234} = 1$ can not be written as $x_1a_1$ with these nonzero components $a_{1,14}$. As a fact, we will consider the only case $a_{1,13} = 0$ in our evaluation below. Start to compute all products $a_ia_j$, $i, j = 2, 3, \ldots 15$.

$a_2a_2 = (e_1 + a_{2,13}e_{124})(e_1 + a_{2,13}e_{124}) = -1 - a_{2,13}e_{24} - a_{2,13}e_{24} + a_{2,13}^2 = x_1a_1 + x_{10}a_{10}$. So, $x_1 = a_{2,13}^2 - 1$, $x_{10} = -2a_{2,13}$, and $2a_{2,13}a_{10,13} = 0$. It means that $a_{2,13} = 0$ or $a_{10,13} = 0$.

$a_2a_3 = (e_1 + a_{2,13}e_{124})(e_2 + a_{3,13}e_{124}) = e_{12} - a_{3,13}e_{24} + a_{2,13}e_{14} + a_{2,13}a_{3,13} = x_1a_1 + x_6a_6 + x_8a_8 + x_{10}a_{10}$. So, $x_1 = a_{2,13}a_{3,13}$, $x_6 = 1$, $x_8 = a_{2,13}$, $x_{10} = -a_{3,13}$, and $a_{6,13} + a_{2,13}a_{8,13} - a_{3,13}a_{10,13} = 0$. Vice versa, $a_3a_2 = (e_2 + a_{3,13}e_{124})(e_1 + a_{2,13}e_{124}) = -e_{12} + a_{2,13}e_{14} - a_{3,13}e_{24} + a_{2,13}a_{3,13} = x_1a_1 + x_6a_6 + x_8a_8 + x_{10}a_{10}$. So, $x_1 = a_{2,13}a_{3,13}$, $x_6 = -1$, $x_8 = a_{2,13}$, $x_{10} = -a_{3,13}$, and $-a_{6,13} + a_{2,13}a_{8,13} - a_{3,13}a_{10,13} = 0$. Last two equations produce the following results $a_{6,13} = 0$, $a_{2,13}a_{8,13} - a_{3,13}a_{10,13} = 0$.

$a_2a_4 = (e_1 + a_{2,13}e_{124})(e_3 + a_{4,13}e_{124}) = e_{13} - a_{4,13}e_{24} - a_{2,13}e_{1234} + a_{2,13}a_{4,13} = x_1a_1 + x_7a_7 + x_{10}a_{10} + x_{15}a_{15}$. So, $x_1 = a_{2,13}a_{4,13}$, $x_7 = 1$, $x_{10} = -a_{4,13}$, $x_{15} = -a_{2,13}$, and $a_{7,13} + a_{4,13}a_{10,13} = 0$. Vice versa, $a_4a_2 = (e_3 + a_{4,13}e_{124})(e_1 + a_{2,13}e_{124}) = -e_{13} - a_{4,13}e_{24} + a_{2,13}e_{1234} + a_{2,13}a_{4,13} = x_1a_1 + x_7a_7 + x_{10}a_{10} + x_{15}a_{15}$. So, $x_1 = a_{2,13}a_{4,13}$, $x_7 = 1$, $x_{10} = -a_{4,13}$, $x_{15} = -a_{2,13}$, and $a_{7,13} - a_{4,13}a_{10,13} = 0$. Last two equations produce the following results $a_{7,13} = 0$, $a_{4,13}a_{10,13} = 0$.

$a_2a_5 = (e_1 + a_{2,13}e_{124})(e_4 + a_{5,13}e_{124}) = e_{14} - a_{5,13}e_{24} - a_{2,13}e_{12} + a_{2,13}a_{5,13} = x_1a_1 + x_6a_6 + x_8a_8 + x_{10}a_{10}$. So, $x_1 = a_{2,13}a_{5,13}$, $x_6 = -a_{2,13}$, $x_8 = 1$, $x_{10} = -a_{5,13}$, and $-a_{2,13}a_{6,13} + a_{8,13} - a_{5,13}a_{10,13} = 0$. Vice versa, $a_5a_2 = (e_4 + a_{5,13}e_{124})(e_1 + a_{2,13}e_{124}) = -e_{14} - a_{5,13}e_{24} - a_{2,13}e_{12} + a_{2,13}a_{5,13} = x_1a_1 + x_6a_6 + x_8a_8 + x_{10}a_{10}$. So, $x_1 = a_{2,13}a_{5,13}$, $x_6 = -a_{2,13}$, $x_8 = -1$, $x_{10} = -a_{5,13}$, and $-a_{2,13}a_{6,13} - a_{8,13} - a_{5,13}a_{10,13} = 0$. Last two equations produce the following results $a_{8,13} = 0$, $a_{2,13}a_{6,13} = -a_{5,13}a_{10,13}$.

$a_2a_6 = (e_1 + a_{2,13}e_{124})(e_{12} + a_{6,13}e_{124}) = -e_2 - a_{6,13}e_{24} - a_{2,13}e_4 + a_{2,13}a_{6,13} = x_1a_1 + x_3a_3 + x_5a_5 + x_{10}a_{10}$. So, $x_1 = a_{2,13}a_{6,13}$, $x_3 = -1$, $x_5 = -a_{2,13}$, $x_{10} = -a_{6,13}$, and $-a_{3,13} - a_{2,13}a_{5,13} - a_{6,13}a_{10,13} = 0$. Vice versa, $a_6a_2 = (e_{12} + a_{6,13}e_{124})(e_1 + a_{2,13}e_{124}) = e_2 - a_{6,13}e_{24} - a_{2,13}e_4 + a_{2,13}a_{6,13} = x_1a_1 + x_3a_3 + x_5a_5 + x_{10}a_{10}$. So, $x_1 = a_{2,13}a_{6,13}$, $x_3 = 1$, $x_5 = -a_{2,13}$, $x_{10} = -a_{6,13}$, and $a_{3,13} - a_{2,13}a_{5,13} - a_{6,13}a_{10,13} = 0$. Last two equations produce the following results $a_{3,13} = 0$, $a_{2,13}a_{5,13} = -a_{6,13}a_{10,13}$

$a_2a_7 = (e_1 + a_{2,13}e_{124})(e_{13} + a_{7,13}e_{124}) = -e_3 - a_{7,13}e_{24} + a_{2,13}e_{234} + a_{2,13}a_{7,13} = x_1a_1 + x_4a_4 + x_{10}a_{10} + x_{14}a_{14}$. So, $x_1 = a_{2,13}a_{7,13}$, $x_4 = -1$, $x_{10} = -a_{7,13}$, $x_{14} = a_{2,13}$, and $-a_{4,13} - a_{7,13}a_{10,13} = 0$. Vice versa, $a_7a_2 = (e_{13} + a_{7,13}e_{124})(e_1 + a_{2,13}e_{124}) = e_3 - a_{7,13}e_{24} - a_{2,13}e_{234} + a_{2,13}a_{7,13} = x_1a_1 + x_4a_4 + x_{10}a_{10} + x_{14}a_{14}$. So, $x_1 = a_{2,13}a_{7,13}$, $x_4 = 1$, $x_{10} = -a_{7,13}$, $x_{14} = -a_{2,13}$, and $-a_{4,13} - a_{7,13}a_{10,13} = 0$. Last two equations produce the following results $a_{4,13} = 0$, $a_{7,13}a_{10,13} = 0$.

$a_2a_8 = (e_1 + a_{2,13}e_{124})(e_{14} + a_{8,13}e_{124}) = -e_4 - a_{8,13}e_{24} + a_{2,13}e_2 + a_{2,13}a_{8,13} = x_1a_1 + x_3a_3 + x_5a_5 + x_{10}a_{10}$. So, $x_1 = a_{2,13}a_{8,13}$, $x_3 = a_{2,13}$, $x_5 = -1$, $x_{10} = -a_{8,13}$, and $a_{2,13}a_{3,13} - a_{5,13} - a_{8,13}a_{10,13} = 0$. Product $a_8a_2$ gives: $a_8a_2 = (e_{14} + a_{8,13}e_{124})(e_1 +$

$a_{2,13}e_{124}) = e_4 + a_{2,13}e_2 - a_{8,13}e_{24} + a_{2,13}a_{8,13} = x_1a_1 + x_3a_3 + x_5a_5 + x_{10}a_{10}$. So, $x_1 = a_{2,13}a_{8,13}$, $x_3 = a_{2,13}$, $x_5 = 1$, $x_{10} = -a_{8,13}$, and $a_{2,13}a_{3,13} + a_{5,13} - a_{8,13}a_{10,13} = 0$.. Last two equations produce the following results $a_{5,13} = 0$, $a_{2,13}a_{3,13} = a_{8,13}a_{10,13}$.

$a_2a_9 = (e_1 + a_{2,13}e_{124})(e_{23} + a_{9,13}e_{124}) = e_{123} - a_{9,13}e_{24} - a_{2,13}e_{134} + a_{2,13}a_{9,13} = x_1a_1 + x_{10}a_{10} + x_{12}a_{12} + x_{13}a_{13}$. So, $x_1 = a_{2,13}a_{9,13}$, $x_{10} = -a_{9,13}$, $x_{12} = 1$, $x_{13} = -a_{2,13}$, and $-a_{9,13}a_{10,13} + a_{12,13} = 0$. Vice versa, $a_9a_2 = (e_{23} + a_{9,13}e_{124})(e_1 + a_{2,13}e_{124}) = e_{123} - a_{9,13}e_{24} + a_{2,13}e_{134} + a_{2,13}a_{9,13} = x_1a_1 + x_{10}a_{10} + x_{12}a_{12} + x_{13}a_{13}$. So, $x_1 = a_{2,13}a_{9,13}$, $x_{10} = -a_{9,13}$, $x_{12} = 1$, $x_{13} = a_{2,13}$, and $-a_{9,13}a_{10,13} + a_{12,13} = 0$. Last two equations coinside, we have the same one equation $a_{12,13} = a_{9,13}a_{10,13}$.

$a_2a_{10} = (e_1 + a_{2,13}e_{124})(e_{24} + a_{10,13}e_{124}) = e_{124} - a_{10,13}e_{24} - a_{2,13}e_1 + a_{2,13}a_{10,13} = x_1a_1 + x_2a_2 + x_{10}a_{10}$. So, $x_1 = a_{2,13}a_{10,13}$, $x_2 = -a_{2,13}$, $x_{10} = -a_{10,13}$, and $-a_{10,13}a_{10,13} - a_{2,13}a_{2,13} = 1$. Vice versa, $a_{10}a_2 = (e_{24} + a_{10,13}e_{124})(e_1 + a_{2,13}e_{124}) = e_{124} - a_{10,13}e_{24} - a_{2,13}e_1 + a_{2,13}a_{10,13} = x_1a_1 + x_2a_2 + x_{10}a_{10}$. So, $x_1 = a_{2,13}a_{10,13}$, $x_2 = -a_{2,13}$, $x_{10} = -a_{10,13}$, and we obtain the same equation $-a_{10,13}a_{10,13} - a_{2,13}a_{2,13} = 1$.

$a_2a_{11} = (e_1 + a_{2,13}e_{124})(e_{34} + a_{11,13}e_{124}) = e_{134} - a_{11,13}e_{24} + a_{2,13}e_{123} + a_{2,13}a_{11,13} = x_1a_1 + x_{10}a_{10} + x_{12}a_{12} + x_{13}a_{13}$. So, $x_1 = a_{2,13}a_{11,13}$, $x_{10} = -a_{11,13}$, $x_{12} = a_{2,13}$, $x_{13} = 1$, and $-a_{11,13}a_{10,13} + a_{2,13}a_{12,13} = 0$. Vice versa, $a_{11}a_2 = (e_{34} + a_{11,13}e_{124})(e_1 + a_{2,13}e_{124}) = e_{134} - a_{11,13}e_{24} - a_{2,13}e_{123} + a_{2,13}a_{11,13} = x_1a_1 + x_{10}a_{10} + x_{12}a_{12} + x_{13}a_{13}$. . So, $x_1 = a_{2,13}a_{11,13}$, $x_{10} = -a_{11,13}$, $x_{12} = -a_{2,13}$, $x_{13} = 1$, and $-a_{11,13}a_{10,13} + a_{2,13}a_{12,13} = 0$. Last two equations produce the following results $a_{2,13}a_{12,13} = 0$, $a_{11,13}a_{10,13} = 0$.

$a_2a_{12} = (e_1 + a_{2,13}e_{124})(e_{123} + a_{12,13}e_{124}) = -e_{23} - a_{12,13}e_{24} + a_{2,13}e_{34} + a_{2,13}a_{12,13} = x_1a_1 + x_9a_9 + x_{10}a_{10} + x_{11}a_{11}$. So, $x_1 = a_{2,13}a_{12,13}$, $x_9 = -1$, $x_{10} = -a_{12,13}$, $x_{11} = a_{2,13}$, and $-a_{9,13} - a_{12,13}a_{10,13} + a_{2,13}a_{11,13} = 0$. Compute $a_{12}a_2 = (e_{123} + a_{12,13}e_{124})(e_1 + a_{2,13}e_{124}) = -e_{23} - a_{12,13}e_{24} - a_{2,13}e_{34} + a_{2,13}a_{12,13} = x_1a_1 + x_9a_9 + x_{10}a_{10} + x_{11}a_{11}$. So, $x_1 = a_{2,13}a_{12,13}$, $x_9 = -1$, $x_{10} = -a_{12,13}$, $x_{11} = -a_{2,13}$, and $-a_{9,13} - a_{12,13}a_{10,13} - a_{2,13}a_{11,13} = 0$. Last two equations produce the following results $a_{9,13} = -a_{12,13}a_{10,13}$, $a_{2,13}a_{10,13} = 0$.

$a_2a_{13} = (e_1 + a_{2,13}e_{124})(e_{134}) = -e_{34} - a_{2,13}e_{23} = x_9a_9 + x_{11}a_{11}$. So, $x_9 = -a_{2,13}$, $x_{11} = -1$, and $-a_{2,13}a_{9,13} - a_{11,13} = 0$. Compute now $a_{13}a_2 = (e_{134})(e_1 + a_{2,13}e_{124}) = -e_{34} + a_{2,13}e_{23} = x_9a_9 + x_{11}a_{11}$. . So, $x_9 = a_{2,13}$, $x_{11} = -1$, and $a_{2,13}a_{9,13} - a_{11,13} = 0$. Last two equations produce the following results $a_{11,13} = 0$, $a_{2,13}a_{9,13} = 0$.

$a_2a_{14} = (e_1 + a_{2,13}e_{124})(e_{234}) = e_{1234} + a_{2,13}e_{13} = x_7a_7 + x_{15}a_{15}$. So, $x_7 = a_{2,13}$, $x_{15} = 1$, and $a_{2,13}a_{7,34} = 0$. Product $a_{14}a_2$ generates the same equation.

$a_2a_{15} = (e_1 + a_{2,13}e_{124})(e_{1234}) = -e_{234} - a_{2,13}e_3 = x_4a_4 + x_{14}a_{14}$. So, $x_4 = -a_{2,13}$, $x_{14} = -1$, and $a_{2,13}a_{4,13} = 0$.

Our evaluation shows that $a_{3,13} = 0$, $a_{4,13} = 0$, $a_{5,13} = 0$, $a_{6,13} = 0$, $a_{7,13} = 0$, $a_{8,13} = 0$, and $a_{11,13} = 0$ in particular. Compute now the following product $a_5a_6 = e_4e_{12} = e_{124}$. But this element $e_{124}$ can not be written as a linear combination of vectors from the Basis (4). It means that basis (4) doesn't generate any 15-dimensional subalgebra in this 16-dimensional Clifford algebra.

**Basis (5):** $a_1 = 1 + a_{1,12}e_{123}$, $a_2 = e_1 + a_{2,12}e_{123}$, $a_3 = e_2 + a_{3,12}e_{123}$, $a_4 = e_3 + a_{4,12}e_{123}$, $a_5 = e_4 + a_{5,12}e_{123}$, $a_6 = e_{12} + a_{6,12}e_{123}$, $a_7 = e_{13} + a_{7,12}e_{123}$, $a_8 = e_{14} + a_{8,12}e_{123}$, $a_9 = e_{23} + a_{9,12}e_{123}$, $a_{10} = e_{24} + a_{10,12}e_{123}$, $a_{11} = e_{34} + a_{11,12}e_{123}$, $a_{12} = e_{124}$, $a_{13} = e_{134}$, $a_{14} = e_{234}$, $a_{15} = e_{1234}$.

First of all, compute $a_1 a_1 = \left(1 + a_{1,12}e_{123}\right)\left(1 + a_{1,12}e_{123}\right) = 1 + a_{1,12}^2 + 2a_{1,12}e_{123} = x_1 a_1$. So, $x_1 = 1 + a_{1,12}^2$, and $a_{1,12}\left(a_{1,12}^2 - 1\right) = 0$. Therefore, $a_{1,12} = 0$, or $a_{1,12} = 1$, or $a_{1,12} = -1$. Last two cases, $a_{1,12} = 1$ and $a_{1,12} = -1$, are impossible because the product $a_{15}a_{15} = e_{1234}e_{1234} = 1$ can not be written as $x_1 a_1$ with those nonzero components $a_{1,12}$. As a fact, we will consider the only case $a_{1,12} = 0$ in our evaluation below. Start to compute all products $a_i a_j$, $i, j = 2, 3, \ldots 15$.

$a_2 a_2 = (e_1 + a_{2,12}e_{123})(e_1 + a_{2,12}e_{123}) = -1 - a_{2,12}e_{23} - a_{2,12}e_{23} + a_{2,13}^2 = x_1 a_1 + x_9 9$. So, $x_1 = a_{2,12}^2 - 1$, $x_9 = -2a_{2,12}$, and $2a_{2,12}a_{9,12} = 0$. It means that $a_{2,12} = 0$ or $a_{9,12} = 0$.

$a_2 a_3 = (e_1 + a_{2,12}e_{123})(e_2 + a_{3,12}e_{123}) = e_{12} - a_{3,12}e_{23} + a_{2,12}e_{13} + a_{2,12}a_{3,12} = x_1 a_1 + x_6 a_6 + x_7 a_7 + x_9 a_9$. So, $x_1 = a_{2,12}a_{3,12}$, $x_6 = 1$, $x_7 = a_{2,12}$, $x_9 = -a_{3,12}$, and $a_{6,12} + a_{2,12}a_{7,12} - a_{3,12}a_{9,12} = 0$. Vice versa, $a_3 a_2 = (e_2 + a_{3,12}e_{123})(e_1 + a_{2,12}e_{123}) = -e_{12} + a_{2,12}e_{13} - a_{3,12}e_{23} + a_{2,12}a_{3,12} = x_1 a_1 + x_6 a_6 + x_7 a_7 + x_9 a_9$. So, $x_1 = a_{2,12}a_{3,12}$, $x_6 = -1$, $x_7 = a_{2,12}$, $x_9 = -a_{3,12}$, and $-a_{6,12} + a_{2,12}a_{7,12} - a_{3,12}a_{9,12} = 0$. Last two equations produce the following results $a_{6,12} = 0$, $a_{2,12}a_{7,12} - a_{3,12}a_{9,12} = 0$.

$a_2 a_4 = (e_1 + a_{2,12}e_{123})(e_3 + a_{4,12}e_{123}) = e_{13} - a_{4,12}e_{23} - a_{2,12}e_{12} + a_{2,12}a_{4,12} = x_1 a_1 + x_6 a_6 + x_7 a_7 + x_9 a_9$. So, $x_1 = a_{2,12}a_{4,12}$, $x_6 = -a_{2,12}$, $x_7 = 1$, $x_9 = -a_{4,12}$, and $-a_{2,12}a_{6,12} + a_{7,12} + a_{4,12}a_{9,12} = 0$. Vice versa, $a_4 a_2 = (e_3 + a_{4,12}e_{123})(e_1 + a_{2,12}e_{123}) = -e_{13} - a_{4,12}e_{23} - a_{2,12}e_{12} + a_{2,12}a_{4,12} = x_1 a_1 + x_6 a_6 + x_7 a_7 + x_9 a_9$. So, $x_1 = a_{2,12}a_{4,12}$, $x_6 = -a_{2,12}$, $x_7 = -1$, $x_9 = -a_{4,12}$, and $-a_{2,12}a_{6,12} - a_{7,12} + a_{4,12}a_{9,12} = 0$. Last two equations produce the following results $a_{7,12} = 0$, $a_{2,12}a_{6,12} = a_{4,12}a_{9,12}$.

$a_2 a_5 = (e_1 + a_{2,12}e_{123})(e_4 + a_{5,12}e_{123}) = e_{14} - a_{5,12}e_{23} + a_{2,12}e_{1234} + a_{2,12}a_{5,12} = x_1 a_1 + x_8 a_8 + x_9 a_9 + x_{15}a_{15}$. So, $x_1 = a_{2,12}a_{5,12}$, $x_8 = 1$, $x_9 = -a_{5,12}$, $x_{15} = a_{2,12}$, and $a_{8,12} - a_{5,12}a_{9,12} = 0$. Vice versa, $a_5 a_2 = (e_4 + a_{5,12}e_{123})(e_1 + a_{2,12}e_{123}) = -e_{14} - a_{5,12}e_{23} - a_{2,12}e_{1234} + a_{2,12}a_{5,12} = x_1 a_1 + x_8 a_8 + x_9 a_9 + x_{15}a_{15}$. So, $x_1 = a_{2,12}a_{5,12}$, $x_8 = -1$, $x_9 = -a_{5,12}$, $x_{15} = -a_{2,12}$, and $-a_{8,12} - a_{5,12}a_{9,12} = 0$. Last two equations produce the following results $a_{8,12} = 0$, $a_{5,12}a_{9,12} = 0$.

$a_2 a_6 = (e_1 + a_{2,12}e_{123})(e_{12} + a_{6,12}e_{123}) = -e_2 - a_{6,12}e_{23} - a_{2,12}e_3 + a_{2,12}a_{6,12} = x_1 a_1 + x_3 a_3 + x_4 a_4 + x_9 a_9$. So, $x_1 = a_{2,12}a_{6,12}$, $x_3 = -1$, $x_4 = -a_{2,12}$, $x_9 = -a_{6,12}$, and $-a_{3,12} - a_{2,13}a_{4,12} - a_{6,12}a_{9,12} = 0$. Vice versa, $a_6 a_2 = (e_{12} + a_{6,12}e_{123})(e_1 + a_{2,12}e_{123}) = e_2 - a_{6,12}e_{23} - a_{2,12}e_3 + a_{2,12}a_{6,12} = x_1 a_1 + x_3 a_3 + x_4 a_4 + x_9 a_9$. So, $x_1 = a_{2,12}a_{6,12}$, $x_3 = 1$, $x_4 = -a_{2,12}$, $x_9 = -a_{6,12}$, and $a_{3,12} - a_{2,13}a_{4,12} - a_{6,12}a_{9,12} = 0$. Last two equations produce the following results $a_{3,12} = 0$, $a_{2,13}a_{4,12} = -a_{6,12}a_{9,12}$.

$a_2 a_7 = (e_1 + a_{2,12}e_{123})(e_{13} + a_{7,12}e_{123}) = -e_3 - a_{7,12}e_{23} + a_{2,12}e_2 + a_{2,12}a_{7,12} = x_1 a_1 + x_3 a_3 + x_4 a_4 + x_9 a_9$. So, $x_1 = a_{2,12}a_{7,12}$, $x_3 = a_{2,12}$, $x_4 = -1$, $x_9 = -a_{7,12}$, and $a_{2,12}a_{3,12} - a_{4,12} - a_{7,12}a_{9,12} = 0$. Vice versa, $a_7 a_2 = (e_{13} + a_{7,12}e_{123})(e_1 + a_{2,12}e_{123}) = e_3 -$

$a_{7,12}e_{23} + a_{2,12}e_2 + a_{2,12}a_{7,12} = x_1a_1 + x_3a_3 + x_4a_4 + x_9a_9$. So, $x_1 = a_{2,12}a_{7,12}$, $x_3 = a_{2,12}$, $x_4 = 1$, $x_9 = -a_{7,12}$, and $a_{2,12}a_{3,12} + a_{4,12} - a_{7,12}a_{9,12} = 0$. Last two equations produce the following results $a_{4,12} = 0$, $a_{2,12}a_{3,12} = a_{7,12}a_{9,12}$.

$a_2a_8 = (e_1 + a_{2,12}e_{123})(e_{14} + a_{8,12}e_{123}) = -e_4 - a_{8,12}e_{23} + a_{2,12}e_{234} + a_{2,12}a_{8,12} = x_1a_1 + x_5a_5 + x_9a_9 + x_{14}a_{14}$. So, $x_1 = a_{2,12}a_{8,12}$, $x_5 = -1$, $x_9 = -a_{8,12}$, $x_{14} = a_{2,12}$, and $-a_{5,12} - a_{8,12}a_{9,12} = 0$. Product $a_8a_2$ gives: $a_8a_2 = (e_{14} + a_{8,12}e_{123})(e_1 + a_{2,12}e_{123}) = e_4 + a_{2,12}e_{234} - a_{8,12}e_{23} + a_{2,12}a_{8,12} = x_1a_1 + x_5a_5 + x_9a_9 + x_{14}a_{14}$. So, $x_1 = a_{2,12}a_{8,12}$, $x_5 = 1$, $x_9 = -a_{8,12}$, $x_{14} = a_{2,12}$, and $-a_{5,12} - a_{8,12}a_{9,12} = 0$. Last two equations produce the following results $a_{5,12} = 0$, $a_{8,12}a_{9,12} = 0$.

$a_2a_9 = (e_1 + a_{2,12}e_{123})(e_{23} + a_{9,12}e_{123}) = e_{123} - a_{9,12}e_{23} - a_{2,12}e_1 + a_{2,12}a_{9,12} = x_1a_1 + x_2a_2 + x_9a_9$. So, $x_1 = a_{2,12}a_{9,12}$, $x_2 = -a_{2,12}$, $x_9 = -a_{9,12}$, and $-a_{2,12}^2 - a_{9,12}^2 = 1$. Vice versa, $a_9a_2 = (e_{23} + a_{9,12}e_{123})(e_1 + a_{2,12}e_{123}) = e_{123} - a_{9,12}e_{23} - a_{2,12}e_1 + a_{2,12}a_{9,12} = x_1a_1 + x_2a_2 + x_9a_9$. So, $x_1 = a_{2,12}a_{9,12}$, $x_2 = -a_{2,12}$, $x_9 = -a_{9,12}$, and $-a_{2,12}^2 - a_{9,12}^2 = 1$. The same one equation is received.

$a_2a_{10} = (e_1 + a_{2,12}e_{123})(e_{24} + a_{10,12}e_{123}) = e_{124} - a_{10,12}e_{23} + a_{2,12}e_{134} + a_{2,12}a_{10,12} = x_1a_1 + x_9a_9 + x_{12}a_{12} + x_{13}a_{13}$. So, $x_1 = a_{2,12}a_{10,12}$, $x_9 = -a_{10,12}$, $x_{12} = 1$, $x_{13} = a_{2,12}$, and $-a_{9,12}a_{10,12} = 0$. Vice versa, $a_{10}a_2 = (e_{24} + a_{10,12}e_{123})(e_1 + a_{2,12}e_{123}) = e_{124} - a_{10,12}e_{23} - a_{2,12}e_{134} + a_{2,12}a_{10,12} = x_1a_1 + x_9a_9 + x_{12}a_{12} + x_{13}a_{13}$. So, $x_1 = a_{2,12}a_{10,12}$, $x_9 = -a_{10,12}$, $x_{12} = 1$, $x_{13} = -a_{2,12}$, and $-a_{9,12}a_{10,12} = 0$. The same one equation is exactly received.

$a_2a_{11} = (e_1 + a_{2,12}e_{123})(e_{34} + a_{11,12}e_{123}) = e_{134} - a_{11,12}e_{23} - a_{2,12}e_{124} + a_{2,12}a_{11,12} = x_1a_1 + x_9a_9 + x_{12}a_{12} + x_{13}a_{13}$. So, $x_1 = a_{2,12}a_{11,12}$, $x_9 = -a_{11,12}$, $x_{12} = -a_{2,12}$, $x_{13} = 1$, and $-a_{11,12}a_{9,12} = 0$. Vice versa, $a_{11}a_2 = (e_{34} + a_{11,12}e_{123})(e_1 + a_{2,12}e_{123}) = e_{134} - a_{11,12}e_{23} + a_{2,12}e_{124} + a_{2,12}a_{11,12} = x_1a_1 + x_9a_9 + x_{12}a_{12} + x_{13}a_{13}$. . So, $x_1 = a_{2,12}a_{11,12}$, $x_9 = -a_{11,12}$, $x_{12} = a_{2,12}$, $x_{13} = 1$, and $-a_{9,12}a_{11,12} = 0$. The same one equation is exactly received.

$a_2a_{12} = (e_1 + a_{2,12}e_{123})(e_{124}) = -e_{24} - a_{2,12}e_{34} = x_{10}a_{10} + x_{11}a_{11}$. So, $x_{10} = -1$, $x_{11} = -a_{2,12}$, and $-a_{10,12} - a_{2,12}a_{11,12} = 0$. Compute $a_{12}a_2 = (e_{124})(e_1 + a_{2,12}e_{123}) = -e_{24} + a_{2,12}e_{34} = x_{10}a_{10} + x_{11}a_{11}$. So, $x_{10} = -1$, $x_{11} = a_{2,12}$, and $-a_{10,12} + a_{2,12}a_{11,12} = 0$. Last two equations produce the following results $a_{10,12} = 0$, $a_{2,12}a_{11,12} = 0$.

$a_2a_{13} = (e_1 + a_{2,12}e_{123})(e_{134}) = -e_{34} + a_{2,12}e_{24} = x_{10}a_{10} + x_{11}a_{11}$. So, $x_{10} = a_{2,12}$, $x_{11} = -1$, and $a_{2,12}a_{10,12} - a_{11,12} = 0$. Compute now $a_{13}a_2 = (e_{134})(e_1 + a_{2,12}e_{123}) = -e_{34} - a_{2,12}e_{24} = x_{10}a_{10} + x_{11}a_{11}$. . So, $x_{10} = -a_{2,12}$, $x_{11} = -1$, and $-a_{2,12}a_{10,12} - a_{11,12} = 0$. Last two equations produce the following results $a_{11,12} = 0$, $a_{2,12}a_{10,12} = 0$.

$a_2a_{14} = (e_1 + a_{2,12}e_{123})(e_{234}) = e_{1234} - a_{2,12}e_{14} = x_8a_8 + x_{15}a_{15}$. So, $x_8 = -a_{2,12}$, $x_{15} = 1$, and $a_{2,12}a_{8,12} = 0$. Product $a_{14}a_2$ generates the same equation.

$a_2a_{15} = (e_1 + a_{2,12}e_{123})(e_{1234}) = -e_{234} + a_{2,12}e_4 = x_5a_5 + x_{14}a_{14}$. So, $x_4 = a_{2,12}$, $x_{14} = -1$, and $a_{2,12}a_{5,12} = 0$.

Our evaluation shows that $a_{3,12} = 0$, $a_{4,12} = 0$, $a_{5,12} = 0$, $a_{6,12} = 0$, $a_{7,12} = 0$, $a_{8,12} = 0$, $a_{10,12} = 0$, and $a_{11,12} = 0$. Compute now the following product $a_4 a_6 = e_3 e_{12} = e_{123}$. But this element $e_{123}$ can not be written as a linear combination of vectors from the Basis (5). It means that basis (5) doesn't generate any 15-dimensional subalgebra at this 16-dimensional Clifford algebra.

**Basis (6):** $a_1 = 1 + a_{1,11}e_{34}$, $a_2 = e_1 + a_{2,11}e_{34}$, $a_3 = e_2 + a_{3,11}e_{34}$, $a_4 = e_3 + a_{4,11}e_{34}$, $a_5 = e_4 + a_{5,11}e_{34}$, $a_6 = e_{12} + a_{6,11}e_{34}$, $a_7 = e_{13} + a_{7,11}e_{34}$, $a_8 = e_{14} + a_{8,11}e_{34}$, $a_9 = e_{23} + a_{9,11}e_{34}$, $a_{10} = e_{24} + a_{10,11}e_{34}$, $a_{11} = e_{123}$, $a_{12} = e_{124}$, $a_{13} = e_{134}$, $a_{14} = e_{234}$, $a_{15} = e_{1234}$.

At least one product is very specific for Basis (6): $a_{11}a_{12} = e_{123}e_{124} = -e_{34}$. This element $-e_{34}$ can not be written as a linear combination of vectors $a_1, \ldots, a_{15}$ from Basis (6). It means that Basis (6) doesn't generate any 15-dimensional subalgebra in 16-dimensional Clifford algebra.

**Basis (7):** $a_1 = 1 + a_{1,10}e_{24}$, $a_2 = e_1 + a_{2,10}e_{24}$, $a_3 = e_2 + a_{3,10}e_{24}$, $a_4 = e_3 + a_{4,10}e_{24}$, $a_5 = e_4 + a_{5,10}e_{24}$, $a_6 = e_{12} + a_{6,10}e_{24}$, $a_7 = e_{13} + a_{7,10}e_{24}$, $a_8 = e_{14} + a_{8,10}e_{24}$, $a_9 = e_{23} + a_{9,10}e_{24}$, $a_{10} = e_{34}$, $a_{11} = e_{123}$, $a_{12} = e_{124}$, $a_{13} = e_{134}$, $a_{14} = e_{234}$, $a_{15} = e_{1234}$.

Consider the product $a_{11}a_{13} = e_{123}e_{134} = e_{24}$. Unfortunately, the resulting element $e_{24}$ can not be written as a linear combination of elements $a_1, \ldots, a_{15}$ from Basis (7). It means that Basis (7) doesn't generate any 15-dimensional subalgebra in the given 16-dimensional Clifford algebra.

**Basis (8):** $a_1 = 1 + a_{1,9}e_{23}$, $a_2 = e_1 + a_{2,9}e_{23}$, $a_3 = e_2 + a_{3,9}e_{23}$, $a_4 = e_3 + a_{4,9}e_{23}$, $a_5 = e_4 + a_{5,9}e_{23}$, $a_6 = e_{12} + a_{6,9}e_{23}$, $a_7 = e_{13} + a_{7,9}e_{23}$, $a_8 = e_{14} + a_{8,9}e_{23}$, $a_9 = e_{24}$, $a_{10} = e_{34}$, $a_{11} = e_{123}$, $a_{12} = e_{124}$, $a_{13} = e_{134}$, $a_{14} = e_{234}$, $a_{15} = e_{1234}$.

Consider the product $a_{13}a_{12} = e_{134}e_{124} = e_{23}$. Unfortunately, the resulting element $e_{23}$ can not be written as a linear combination of elements $a_1, \ldots, a_{15}$ from Basis (8). It means that Basis (8) doesn't generate any 15-dimensional subalgebra in the given 16-dimensional Clifford algebra.

**Basis (9):** $a_1 = 1 + a_{1,8}e_{14}$, $a_2 = e_1 + a_{2,8}e_{14}$, $a_3 = e_2 + a_{3,8}e_{14}$, $a_4 = e_3 + a_{4,8}e_{14}$, $a_5 = e_4 + a_{5,8}e_{14}$, $a_6 = e_{12} + a_{6,8}e_{14}$, $a_7 = e_{13} + a_{7,8}e_{14}$, $a_8 = e_{23}$, $a_9 = e_{24}$, $a_{10} = e_{34}$, $a_{11} = e_{123}$, $a_{12} = e_{124}$, $a_{13} = e_{134}$, $a_{14} = e_{234}$, $a_{15} = e_{1234}$.

Consider the product $a_{14}a_{11} = e_{234}e_{123} = e_{14}$. Unfortunately, the resulting element $e_{14}$ can not be written as a linear combination of elements $a_1, \ldots, a_{15}$ from Basis (9). It means that Basis (9) doesn't generate any 15-dimensional subalgebra in the given 16-dimensional Clifford algebra.

**Basis (10):** $a_1 = 1 + a_{1,7}e_{13}$, $a_2 = e_1 + a_{2,7}e_{13}$, $a_3 = e_2 + a_{3,7}e_{13}$, $a_4 = e_3 + a_{4,7}e_{13}$, $a_5 = e_4 + a_{5,7}e_{13}$, $a_6 = e_{12} + a_{6,7}e_{13}$, $a_7 = e_{14}$, $a_8 = e_{23}$, $a_9 = e_{24}$, $a_{10} = e_{34}$, $a_{11} = e_{123}$, $a_{12} = e_{124}$, $a_{13} = e_{134}$, $a_{14} = e_{234}$, $a_{15} = e_{1234}$.

Consider the product $a_{12}a_{14} = e_{124}e_{234} = e_{13}$. Unfortunately, the resulting element $e_{13}$ can not be written as a linear combination of elements $a_1, \ldots, a_{15}$ from Basis (10). It means that Basis (10) doesn't generate any 15-dimensional subalgebra in the given 16-dimensional Clifford algebra.

It is easy to see that for Bases $(11) - (15)$ our situations are similar to those of Bases $(6) - (10)$, and they don't generate any 15-dimensional subalgebra in the given 16-dimensional Clifford algebra. Consider canonical Basis (16) separately.

**Basis (16):** $a_1 = e_1$, $a_2 = e_2$, $a_3 = e_3$, $a_4 = e_4$, $a_5 = e_{12}$, $a_6 = e_{13}$, $a_7 = e_{14}$, $a_8 = e_{23}$, $a_9 = e_{24}$, $a_{10} = e_{34}$, $a_{11} = e_{123}$, $a_{12} = e_{124}$, $a_{13} = e_{134}$, $a_{14} = e_{234}$, $a_{15} = e_{1234}$.

A lot of products produce element 1 or $-1$ like $a_1a_1 = e_1e_1 = -1$. But element 1 is not located in the subspace generated by Basis (16). It means that Basis (16) doesn't generate any 15-dimensional subalgebra in the given 16-dimensional Clifford algebra.

The proof of the Theorem is finished.